%% file: LOTTERYP.tex
\input fonts

\input articlemacros

\input generalmathmacros

\input arrowmacros


\center [submitted to the Annals of Pure and Applied Logic]

\title The Lottery Preparation

\author Joel David Hamkins\cr
	Kobe University and\cr
	The City University of New York\cr

\abstract. The lottery preparation, a new general kind of Laver preparation, works uniformly with supercompact cardinals, strongly compact cardinals, strong cardinals, measurable cardinals, or what have you. And like the Laver preparation, the lottery preparation makes these cardinals indestructible by various kinds of further forcing.  A supercompact cardinal $\k$, for example, becomes fully indestructible by $\ltk$-directed closed forcing; a strong cardinal $\k$ becomes indestructible by $\ltek$-strategically closed forcing; and a strongly compact cardinal $\k$ becomes indestructible by, among others, the forcing to add a Cohen subset to $\k$, the forcing to shoot a club $C\of\k$ avoiding the measurable cardinals and the forcing to add various long Prikry sequences. The lottery preparation works best when performed after fast function forcing, which adds a new completely general kind of Laver function for any large cardinal, thereby freeing the Laver function concept from the supercompact cardinal context.

\footnote{}{\noindent My research has been supported in part by a grant from the PSC-CUNY Research Foundation, by a CUNY Scholar Incentive Award and by a research fellowship from the Japan Society for the Promotion of Science. I am grateful to my gracious hosts here at Kobe University in Japan, and I would like to thank both Arthur Apter for his insightful suggestions concerning my work with strongly compact cardinals and Tadatoshi Miyamoto for his probing questions on an earlier draft of this paper.}The Laver preparation \cite[Lav78], which spectacularly makes any supercompact cardinal $\k$ indestructible by $\ltk$-directed closed forcing, has long been an indispensible tool and recognized as an important milestone in large cardinal set theory. Is there such a preparation for the other large cardinals? The Laver preparation does not seem to work with strongly compact and other cardinals. While strong cardinals are successfully treated in \cite[GitShl89], the fundamental lifting tools currently available fail outrightly when applied to strongly compact non-supercompact cardinals. The technology has simply not been available to make strongly compact cardinals even partly indestructible. The extent of our ignorance is made strikingly plain by the fact that the following question has remained open:

\quiet\question. Can any strongly compact cardinal $\k$ be made indestructible by the forcing $\add(\k,1)$ which adds, by initial segments, a Cohen subset to $\k$?

In this paper I provide a new technology to answer the above question, and to answer it the way that we all hoped it would be answered: any strongly compact cardinal $\k$ can be made indestructible by $\add(\k,1)$ and more. And the technique is limited to neither strongly compact cardinals nor the particular poset $\add(\k,1)$. Specifically, I present here the lottery preparation, a new kind of Laver preparation, which works uniformly with strongly compact cardinals, supercompact cardinals, measurable cardinals, strong cardinals, or what have you, and makes them all indestructible by a variety of forcing notions.

\quiet\theorem Main Lottery Preparation Theorem. The lottery preparation makes a variety of large cardinals indestructible by various forcing notions. Specifically:  
\points 1. The lottery preparation of a supercompact cardinal $\k$ makes the supercompactness of $\k$ indestructible by any $\ltk$-directed closed forcing.\cr
	2. The lottery preparation of a strongly compact cardinal $\k$ makes the strong compactness of $\k$ indestructible by, among others, the forcing $\add(\k,1)$ which adds a Cohen subset to $\k$, the forcing which shoots a club $C\of\k$ avoiding the measurable cardinals and the forcing which adds certain long Prikry sequences.\cr
       3. The lottery preparation of a strong cardinal $\k$ satisfying $2^\k=\k\plus$ makes the strongness of $\k$ indestructible by, among others, any $\ltek$-strategically closed forcing and by $\add(\k,1)$.\cr
        4. With a dash of the \GCH, level-by-level results hold for partially supercompact and partially strong cardinals.\cr

The precise details are in section four. The lottery preparation, which is defined relative to a function $f\from\k\to\k$, works best when the values of $j(f)(\k)$ can be made large for the desired kind of large cardinal embedding. Since fast function forcing adds a generic function $f$ for which the values of $j(f)(\k)$ can be almost arbitrarily specified, the lottery preparation works especially well when performed after fast function forcing and defined relative to this generic fast function. An interesting related result is the fact that fast function forcing adds a new completely general kind of Laver function: 

\quiet\theorem Generalized Laver Function Theorem. After fast function forcing there is a function $\ell\from\k\to(V[f])_\k$ such that for any embedding $j:V[f]\to M[j(f)]$ with critical point $\k$ (whether internal or external) and any $z$ in $H(\l\plus)^{M[j(f)]}$ where $\l=j(f)(\k)$ there is another embedding $j^*:V[f]\to M[j^*(f)]$ such that: 
\points 1. $j^*(\ell)(\k)=z$, \cr
            2. $M[j^*(f)]=M[j(f)]$,\cr
            3. $j^*\restrict V=j\restrict V$, and\cr
            4. If $j$ is the ultrapower by a standard measure $\eta$ concentrating on a set in $V$, then $j^*$ is the ultrapower by a measure $\eta^*$ concentrating on the same set and moreover $\eta^*\intersect V=\eta\intersect V$ and $[\id]_{\eta^*}=[\id]_\eta$.\cr

The standard measures include, among many others, all normal measures, all supercompactness measures and, up to isomorphism, all strong compactness measures. The restriction that $z\in H(\l\plus)$ is not onerous because the value of $j(f)(\k)$ is amazingly mutable, and can be almost arbitrarily specified: for any embedding $j:V[f]\to M[j(f)]$ and any $\a<j(\k)$ there is another embedding $j^*:V[f]\to M[j^*(f)]$ such that $j^*(f)(\k)=\a$, $j^*\restrict V=j\restrict V$ and $M[j^*(f)]\of M[j(f)]$. For these reasons, the generic Laver function $\ell$ can be effectively used with almost any kind of large cardinal embedding much as a Laver function is used with a supercompactness embedding. In this way, fast function forcing frees the Laver function concept from the supercompact cardinal context.

Going back at least to \cite[Men74], where several preservation theorems are proved, set theorists have wondered about the possibility of making strongly compact and other cardinals indestructible by forcing; perhaps the lottery preparation provides an answer.  The larger question, though, of precisely how indestructible these cardinals can be made is still very much open. Probably the lottery preparation provides more indestructibility than I will identify in this paper. It is natural to hope that any strongly compact cardinal can be made fully indestructible, perhaps by the usual sort of reverse Easton preparation, an iteration of closed forcing. The sad fact, however, with which I conclude this paper is that such a preparation is simply impossible.

\quiet\theorem Impossibility Theorem. By preparatory forcing which admits a gap below $\k$ (such as any preparation naively resembling the Laver preparation), if the measurability of $\k$ can be made indestructible by $\ltk$-directed closed forcing, then $\k$ must have been supercompact in the ground model. 

The details for this theorem are in section five.

Let me quickly explain the structure of this paper. First, I introduce Woodin's fast function forcing, showing in section one that it preserves a variety of large cardinals and in section two that it adds a new general kind of Laver function. Next, I introduce the lottery preparation, proving in section three that it preserves a variety of large cardinals and in section four that it makes these cardinals indestructible by various further forcing. Lastly, in section five I prove the Impossibility Theorem. Throughout, I try to use standard notation, and argue freely in ZFC. By $p\from A\to B$, I mean that $p$ is a partial function from $A$ to $B$. And if $p$ is a condition in the poset $\P$, then by $\P\restrict p$ I mean the sub-poset $\set{q\in\P\st q\leq p}$. My focus is almost always on $\k$, the large cardinal at hand, and so invariably, the critical point of whatever embedding I am concerned with will be denoted by $\k$. 

\section Fast Function Forcing

Fast function forcing, invented by W. Hugh Woodin (with an infant form due to Robert Solovay), allows one to add a function $f\from\k\to\k$ such that the value of $j(f)(\k)$ can be almost arbitrarily specified for embeddings $j:V[f]\to M[j(f)]$ in the extension. These functions behave, therefore, in a completely general large cardinal context much like Laver's function does in the supercompact cardinal context. Indeed, in the next section I will prove that with a fast function one can obtain a completely general kind of Laver function in a completely general large cardinal setting. And since the existence of Laver functions in the supercompact cardinal context has proved so indispensible---these functions appear in dozens if not hundreds of articles---the generalized generic Laver functions here may find a broad application. So let's begin with fast function forcing.

The fast function forcing notion $\F$ for the cardinal $\k$ consists of conditions $p\from\k\to\k$ such that $\dom(p)\of\inacc$ has size less than $\k$ and if $\g\in\dom(p)$ then $p\image\g\of\g$ and $\card{p\restrict\g}<\g$. The conditions are ordered by inclusion. The (union of the) generic for this forcing is the {\df fast function} $f\from\k\to\k$, a partial function on $\k$. To emphasize the role of $\k$, I will sometimes denote $\F$ by $\F_\k$. 

By $\F_{\l,\k}$ I mean the version of fast function forcing consisting of conditions with domain in $[\l,\k)$.
It is easy to see, by taking the union of conditions, that $\F_{\l,\k}$ is $\ltel$-directed closed: the only apparant difficulty is the support requirement that $\card{p\restrict\g}<\g$ for $\g\in\dom(p)$; but if $\g>\l$ is inaccessible, then a union of size $\l$ of supports of size less than $\g$ still has size less than $\g$, and so the difficulty is easily addressed. 

\lemma Fast Function Factor Lemma. Below the condition $p=\sing{\<\g,\a>}\in\F$, where $\g$ is inaccessible, $\a$ is an ordinal and $\l$ is the next inaccessible beyond $\g$ and $\a$, the fast function forcing poset factors as $\F\restrict p\iso\F_\g\cross\F_{\l,\k}$. 
\ref\FFFactor

\proof If $q\leq p$ then $q\restrict\g\in\F_\g$ and $\dom(q)$ is disjoint from $(\g,\l)$. Thus, the map $q\mapsto\<q\restrict\g,q\restrict[\l,\k)>$ provides the desired isomorphism.\qed

Thus, if $f\from\k\to\k$ is a fast function on $\k$ and $\g\in\dom(f)$, then $f\restrict\g$ is a fast function on $\g$. 
More generally, the same argument shows that if $f$ is a fast function and $f\image\g\of\g$, then $f\restrict\g\cross f\restrict[\g,\k)$ is generic for the poset $\F_\g\cross\F_{\g,\k}$. Note that if $\g$ is regular, then $\card{\F_\g}\leq\g$. 

\remark Remark on Gap Forcing. For a technical reason which will be made clear later, I will at times want to precede the forcing $\F$ with some small forcing, such as adding a Cohen subset to the least inaccessible cardinal. This kind of forcing is generally benign in the large cardinal context, and I will regard this small forcing as a part of fast function forcing whenever the need arises. The primary reason to do so is that forcing of the form $\P_1*\P_2$, where $\card{\P_1}<\d$ and $\P_2$ is $\lted$-strategically closed in $V^{\P_1}$ is said in \cite[Ham98] and \cite[Ham$\infty$] to admit a {\df gap} at $\d$. The Gap Forcing Theorem of \cite[Ham$\infty$], with a forerunner in \cite[Ham98], asserts that after forcing $V[G]$ which admits a gap at $\d<\k$, any embedding $j:V[G]\to M[j(G)]$ for which $M[j(G)]$ is closed under $\d$-sequences in $V[G]$---and this includes any ultrapower embedding on any set, as well as most strongness extender embeddings---is a lift of an embedding from the ground model. That is, $M\of V$ and $j\restrict V:V\to M$ is an embedding which is definable in $V$. Additionally, if $M[G]$ is $\l$-closed in $V[G]$ for some $\l\geq\d$, then $M$ is $\l$-closed in $V$, and in particular, $j\image\l\in M$; and if $j$ is a $\l$-strongness embedding induced by a natural extender, where $\l$ is either a successor ordinal or has cofinality more than $\d$, then $V_\l\of M$. Thus, the results of \cite[Ham$\infty$] show that gap forcing cannot create new measurable cardinals, strong cardinals, strongly compact cardinals, supercompact cardinals, and so on, with level-by-level versions generally available. In order to appeal to this theorem, therefore, in the context of fast function forcing, I will introduce a very low gap by preceding $\F$ by some very small forcing and hereafter regard this small forcing as a part of fast function forcing, though I will mention it only when I want to apply the Gap Forcing Theorem. \halmos\raise2pt\hbox{\fiverm Remark}\ref\FFFGap

Woodin defined fast function forcing and used it with something below a strong cardinal (in \cite[CumWdn] he and Cummings had an embedding $j:V\to M$ such that $M^\k\of M$ and $j(\k)>\k\plusplus$). His argument, which I give below, works equally well with measurable and supercompact cardinals and, in a modified form, with weakly compact cardinals.  A significant contribution of this paper is that fast function forcing works also with strongly compact cardinals. For presentational clarity, I will present the fast function lifting techniques in the large cardinal order, rather than the temporal order in which the theorems were first proved.

\theorem Fast Function Theorem. Fast function forcing preserves all cardinals and cofinalities and does not disturb the continuum function. Consequently, fast function forcing preserves all inaccessible cardinals.
\ref\FFFCard

\proof Suppose that $\g$ is regular in $V$ but has cofinality $\d<\g$ in $V[f]$. Since $\F$ has size $\k$, we may assume that $\g\leq\k$. There are two cases. First, it may happen that $f\image\d\of\d$. In this case, the forcing factors as $\F_\d\cross\F_{\d,\k}$. The initial forcing $\F_\d$, however, is too small to collapse the cofinality of $\g$ and the tail forcing $\F_{\d,\k}$ is $\lted$-distributive and so cannot collapse the cofinality of $\g$ to $\d$; this contradicts our assumption. Second, alternatively, it may happen that $f(\b)>\d$ for some $\b<\d$. In this case, the forcing factors as $\F_\b\cross\F_{f(\b),\k}$, and again the initial forcing is too small to collapse the cofinality of $\g$ and the tail forcing is too distributive to collapse it to $\d$; so again we reach a contradiction. Thus, fast function forcing preserves all cardinals and cofinalities. 

A similar argument shows that fast function forcing preserves the values of $2^\d$ calculated in the ground model. Again we split into the two cases. If $f\image\d\of\d$, then we may factor the forcing as $\F_\d\cross\F_{\d,\k}$. The initial forcing has size $\d$ and the tail forcing is $\lted$-distributive; so neither can affect the value of $2^\d$. Alternatively, if $f(\b)\geq\d$ for some $\b<\d$, then we may factor the forcing as $\F_\b\cross\F_{f(\b),\k}$ and make the same argument. So the value of $2^\d$ is preserved, and the theorem is proved.\qed

\theorem Fast Function Theorem. Fast function forcing preserves every weakly compact cardinal. Indeed, if $\k$ is weakly compact in $V$, then after adding a fast function $f\from\k\to\k$ there are weakly compact embeddings $j:N[f]\to M[j(f)]$ such that $j(f)(\k)$ is any desired ordinal up to $\k\plus$. 
\ref\FFFWC

\proof Suppose $\k$ is weakly compact in $V$ and $f\from\k\to\k$ is a fast function. Since the previous theorem establishes that $\k$ remains inaccessible in $V[f]$, in order to prove that $\k$ is weakly compact there it suffices only to verify that $\k$ has the tree property. So suppose $T$ is a $\k$-tree in $V[f]$. Thus, $T=\dot T_f$ for some name $\dot T$. In $V$ let $N$ be a transitive elementary substructure of $H(\k\plus)$ of size $\k$ which contains $\F$ and $\dot T$ and is closed under $\ltk$-sequences. Since $\k$ is weakly compact in $V$ there is an embedding $j:N\to M$ with $\cp(j)=\k$. I may assume that $M=\set{j(h)(\k)\st h\in N}$ (by replacing $M$ with the Mostowski collapse of this set if necessary). It follows that $M^{\ltk}\of M$ in $V$ since if $\vec a=\<a_\a\st \a<\b>$ is a sequence in $V$ of elements of $M$ for some $\b<\k$, then build the function $H(\d)=\<h_\a(\d)\st\a<\b>$, where $a_\a=j(h_\a)(\k)$ and $h_\a\in N$. By the closure assumption on $N$ it follows that $H\in N$ and consequently $j(H)\in M$, so $\vec a=\<j(h_\a)(\k)\st\a<\b>\in M$, as desired. Returning to the main argument, now, for any $\a<j(\k)$ the condition $p=\sing{\<\k,\a>}$ is in $j(\F)$. The tail forcing $\F_{\l,j(\k)}$, where $\l$ is the next inaccessible of $M$ beyond $\k$ and $\a$, is $\ltk$-closed in $M$. Since there are only $\k$ many dense sets for this forcing in $M$, we may line them up into a $\k$-sequence and diagonalize to meet them in order to produce in $V$ an $M$-generic $\ftail\of\F_{\l,j(\k)}$. The combination $j(f)=f\union p\union\ftail$ is $M$-generic for $j(\F)$ and consequently the embedding lifts to $j:N[f]\to M[j(f)]$. By construction, $j(f)(\k)=p(\k)=\a$. Now since $T=\dot T_f$ it follows that $T\in N[f]$. Since $T$ is a $\k$-tree it follows that $j(T)$ is a $j(\k)$-tree. Any element on the $\k^\th$ level of $j(T)$ provides a $\k$-branch through $T$. So $\k$ has the tree property in $V[f]$, as desired. And since the orginal $j$ can be chosen so that $j(\k)$ is as large below $\k\plus$ as desired, and $\a$ can be as large below $j(\k)$ as desired, the value of $j(f)(\k)$ can be any ordinal up to $\k\plus$. 

By employing factor arguments as in the previous theorem, it is easy to see more generally that all weakly compact cardinals are preserved.\qed

\theorem Fast Function Theorem.(Woodin) If $2^\k=\k\plus$, then fast function forcing with $\k$ preserves the measurability of $\k$. Indeed, every ultrapower $j:V\to M$ by a measure on $\k$ in $V$ lifts to an ultrapower $j:V[f]\to M[j(f)]$ such that $j(f)(\k)$ is any desired ordinal up to $j(\k)$. Consequently, the lifted embedding can be the ultrapower by a normal measure, even when the original embedding was not. 
\ref\FFMeasurable

\proof Suppose that $\k$ is measurable and $2^\k=\k\plus$ in $V$, that $f\from\k\to\k$ is a fast function and that $j:V\to M$ is the ultrapower embedding by a measure $\mu$ on $\k$ and $\a<j(\k)$. Below the condition $p=\sing{\<\k,\a>}$, factor the forcing as $\F\cross\Ftail$, where $\Ftail=\F_{\l,\k}$ for the next inaccessible $\l$ in $M$ beyond $\k$ and $\a$. The forcing $\Ftail$ is $\ltek$-closed in $M$. Since $2^\k=\k\plus$, a simple counting argument shows that $\card{j(\k\plus)}^V=\k\plus$. Consequently, since $M^\k\of M$, we can line up all the dense sets of $M$ for the forcing $\Ftail$ and diagonalize against them to produce in $V$ an $M$-generic $\ftail$. Thus, in $V[f]$ we may lift the embedding to $j:V[f]\to M[j(f)]$ where $j(f)=f\union p\union\ftail$. So $\k$ remains measurable there, as desired. By construction $j(f)(\k)=p(\k)=\a$.

Some readers may be surprised by the final conclusion of the theorem. Suppose $\a=[\id]_\mu$. By the standard seed techniques (e.g. see \cite[Ham97]), it follows that $M=\set{j(h)(\a)\st h\in V}$. If as above we arrange the lift $j:V[f]\to M[j(f)]$ in such a way that $j(f)(\k)=\a$, then it is easy to see that $M[j(f)]=\set{j(h)(\k)\st h\in V[f]}$, i.e. the seed $\k$ generates the whole embedding. From this, it follows that the lifted embedding $j$ is the ultrapower by the normal measure $\nu=\set{X\st \k\in j(X)}$ in $V[f]$, as desired. I elaborated on this phenomenon in \cite[Ham94].\qed

The factor arguments employed in \FFFCard\ easily extend to show that fast function forcing for $\k$ preserves all measurable cardinals at which the \GCH\ holds. The general phenomenon that the value of $j(f)(\k)$ can be any ordinal up to $j(\k)$ is further explained in the Fast Function Flexibility Theorem below. 

Recall that a cardinal $\k$ is {\df strong} when for every $\l$ it is {\df $\l$-strong}, so that there is an embedding $j:V\to M$ with critical point $\k$ such that $V_\l\of M$. If there is such an embedding, then by factoring through by the canonical extender, there is one such that $M=\set{j(h)(s)\st h\in V\and s\in V_\l}$; one simply replaces $j$ with $\pi\circ j$,  where $\pi$ is the Mostowski collapse of this set. Furthermore, if $\l$ is either a successor ordinal or has cofinality at least $\k$, then for such an embedding $M$ is closed under $\k$-sequences in $V$. 

\theorem Fast Function Theorem.(Woodin) If $2^\k=\k\plus$, then fast function forcing preserves the strongness of $\k$.
\ref\FFFStrong 

\proof The result is completely local, since I will show that if $\k$ is $\l$-strong in $V$ then this is preserved to the fast function extension $V[f]$. Suppose $j:V\to M$ witnesses the $\l$-strongness of $\k$, so that $V_\l\of M$. Let $\d=\card{V_\l}$. Using the canonical extender, I may assume that $M=\set{j(h)(s)\st h\in V\and s\in\d^{<\w}}$. Let $p=\sing{\<\k,\d>}$ be the condition which jumps up to $\d$ at $\k$. Thus, by the Factor Lemma \FFFactor, below $p$ the forcing $j(\F)$ factors as $\F\cross\Ftail$, where $\Ftail$ is $\lted$-closed in $M$. Now use the pair $\<\k,\d>$ as a seed to form the seed hull $\X=\set{j(h)(\k,\d)\st h\in V}\elesub M$ and obtain the factor embedding $$\trianglediagram{V}{j_0}{j}{M_0}{k}{M}$$ where $k:M_0\to M$ is the inverse of the collapse of $\X$. Since $\k$ and $\d$ are in $\X$, it follows that $k(\d_0)=\d$ for some $\d_0<j_0(\k)$, that $k(p_0)=p$ for $p_0=\sing{\<\k,\d_0>}\in j_0(\P)$ and that $\cp(k)>\k$.  The embedding $j_0:V\to M_0$, being generated by the seed $\<\k,\d_0>$, is simply an ultrapower by a measure on $\k$. In particular, since $2^\k=\k\plus$, the diagonalization argument of \FFMeasurable\ provides a lift $j_0:V[f]\to M_0[j_0(f)]$ below the condition $p_0$. It must be that $j_0(f)=f\union p_0\union\ftail^{M_0}$, where $\ftail^{M_0}$ is $M_0$-generic for the ${\lte}\d_0$-closed forcing $\Ftail^{M_0}$. 

I claim that $k$ lifts to $M_0[j_0(f)]$. First, since $cp(k)>\k$, certainly we know that $k$ lifts to $k:M_0[f]\to M[f]$. In order to lift $k$ the rest of the way it suffices to show that $k\image \ftail^{M_0}\of\Ftail$ is $M$-generic. So, suppose $D\in M$ is open and dense in $\Ftail$. Since $M=\set{k(h)(s)\st h\in M_0\and s\in\d^{<\w}}$, it follows that $D=j(\vec D)(s)$ for some $\vec D=\<D_\s\st \s\in\d_0^{<\w}>$ in $M_0$ and $s\in \d^{<\w}$, where every $D_\s$ is an open dense subset of $\Ftail^{M_0}$. Since $\Ftail^{M_0}$ is ${\lte}\d_0$-closed in $M_0$, it follows that $\Dbar =\intersect_\s \vec D_\s\in M_0$ is still open and dense. Furthermore, $k(\Dbar )\of D$. Thus, since $\ftail^{M_0}$ is $M_0$-generic, $k\image \ftail^{M_0}$ meets $D$, as desired. Consequently, $k$ lifts fully to $k:M_0[j_0(f)]\to M[k(j_0(f))]$, where $k(j_0(f))$ is the filter generated by $k\image j_0(f)$. The composition $k\circ j_0$ provides a lift of $j$ to $j:V[f]\to M[j(f)]$. Since $V_\l\of M$ and $f\in M[j(f)]$, it follows that $(V[f])_\l\of M[j(f)]$, and so $\k$ is still $\l$-strong in $V[f]$, as desired.\qed

The next theorem provides the first nontrivial example of the preservation of an arbitrary strongly compact cardinal of which I am aware. I will make a key use of an old technique of Menas \cite[Men74], used also in \cite[Apt98], in order to know that the cardinal remains strongly compact after forcing (Menas and Apter both need a strongly compact limit of supercompact cardinals). Arthur Apter has pointed out that Menas's technique, anachronistically presented in the `dark' ages before the Laver preparation \cite[Lav78], probably had much unrealized potential.  I hope that the results in this paper tend to confirm his view.

\theorem Fast Function Theorem. Fast function forcing preserves the strong compactness of $\k$. Indeed, every strong compactness measure from $V$ extends to a strong compactness measure in the extension.
\ref\FFFSTRC

\proof Suppose that $f\from\k\to\k$ is a $V$-generic fast function, that $\l\geq\k$, and that $\mu_0$ is a fine measure on $P_\k\l$ in $V$. Let $\theta\geq 2^{\l^{<\k}}$, and let $j:V\to M$ be any $\theta$-strongly compact embedding, the ultrapower by a fine measure $\eta$ on $P_\k\theta$ in $V$. By the cover property for strongly compact embeddings, there is a set $Y\in M$ such that $j\image\mu_0\of Y$ and $\card{Y}^M<j(\k)$. I may assume that $Y\of j(\mu_0)$, and consequently $\intersect Y\in j(\mu_0)$. Any element $s_0\in\intersect Y$ is a seed for $\mu_0$ in the sense that $X\in\mu_0\iff s_0\in j(X)$ for any $X\of P_\k\l$ in $V$. Fix such an $s_0$. Let $s=[\id]_\eta$ and $\d=\card{s}^M$, and pick any $\g\geq\d$. Thus, since $\eta$ is a fine measure on $P_\k\theta$, we have $j\image\theta\of s\in j(P_\k\theta)$ and $\theta\leq\d<j(\k)$. Now, in $j(\F)$, let $p$ be the condition $\sing{\<\k,\g>}$. By the Fast Function Factor Lemma, the forcing $j(\F)$ factors below this condition as $\F\cross\Ftail$, where $\Ftail$ is $\lteg$-closed in $M$. Force to add a $V[f]$-generic $\ftail\of\Ftail$, and let $j(f)=f\union p\union\ftail$. By the factorization, this is $M$-generic for $j(\F)$, and consequently the embedding lifts in $V[f][\ftail]$ to $j:V[f]\to M[j(f)]$. By construction, $j(f)(\k)=p(\k)=\g$. Since the forcing $\Ftail$ was $\lteg$-closed in $M$, it is $\lteg$-distributive in $M[f]$; in particular, it adds no new subsets to $\d$ over $M[f]$. Let $\mu_0^*$ be the measure germinated by the seed $s_0$ via the lifted embedding, so that $X\in\mu_0^*\iff s_0\in j(X)$ for $X\in V[f]$. It is easy to see that $\mu_0^*$ measures every set in $V[f]$, that it extends $\mu_0$, and, since $j\image\l\of s_0$, that is fine.  It remains only for me to show that $\mu_0^*\in V[f ]$. For this I will use Menas's key idea in \cite[Men74]. Enumerate in $V$ the nice names $u=\<\dot X_\a\st\a<\theta>$ for the subsets of $P_\k\l$ in $V[f]$ (a simple counting argument shows that there are $2^{\l^{<\k}}$ many of them). Thus, $j(u)\in M$ and consequently also $j(u)\restrict s\in M$. Enumerate $j(u)\restrict s=\<\dot Y_\b\st \b\in s>$, and observe that $\Ydot_{j(\a)}=j(\Xdot_\a)$ for $\a<\theta$. Let
$t=\set{\b\in s\st s_0\in (\Ydot_\b)_{j(f)}}$.  Thus, $t\of s$ and $t\in M[j(f)]$.  Since $s$ has size $\d$ and $\Ftail$ is
$\lted$-distributive in $M[f]$, it follows that $t\in M[f]$, and therefore $t\in V[f]$.  Now simply observe that $(\Xdot_\a)_f\in
\mu_0^*\iff s_0\in j((\Xdot_\a)_f)=j(\Xdot_\a)_{j(f)}=(\Ydot_{j(\a)})_{j(f)}\iff j(\a)\in t$. So $\mu_0^*$ is definable in $V[f]$ from $t$ and $j\restrict\theta$.  Thus, $\mu_0^*\in V[f]$ as desired.\qed

\theorem Fast Function Theorem. Fast function forcing preserves the supercompactness of $\k$; and every supercompactness measure from the ground model extends to a supercompactness measure in the extension.

\proof To see that the supercompactness of $\k$ is preserved, one can simply take $\eta$ to be a $\theta$-supercompactness embedding in the previous argument and $\mu_0$ the $\l$-supercompactness measure germinated via $j$ by the seed $s_0=j\image\l$. The resulting measure $\mu_0^*$ is easily seen to be normal and fine. So $\k$ remains supercompact in $V[f]$.

So now let me show a bit more; namely, that every supercompactness measure from $V$ extends to a supercompactness measure in $V[f]$. Suppose in $V$ that $\mu_0$ is a $\l$-supercompactness measure on $P_\k\l$ and $\eta_0$ is a $\theta$-strong compactness measure on $P_\k\theta$ for some $\theta\geq 2^{\l^\ltk}$. It is not difficult to argue (see the argument preceding Theorem 4.2) that $\eta=\mu_0\cross\eta_0$ is isomorphic to a $\theta$-strong compactness measure whose embedding $j:V\to M$ is closed under $\lambda$-sequences. Furthermore, $s_0=j\image\l$ is a seed for $\mu_0$ via $j$. The previous argument shows how to lift this embedding so that the measure $\mu_0^*$ germinated by $s_0$ via $j:V[f]\to M[j(f)]$ lies in $V[f]$. Again, it is not difficult to argue that $\mu_0^*$ is normal and fine, as desired.\qed

The previous argument actually establishes the following theorem:

\theorem Local Version. If $\k$ is $2^{\l^\ltk}\!$-strongly compact then fast function forcing preserves the $\l$-strong compactness of $\k$. The same holds for supercompactness. Indeed, if $\k$ is $2^{\l^\ltk}\!$-strongly compact and $\l$-supercompact, then fast function forcing preserves the $\l$-supercompactness of $\k$. 
\ref\LocalA

By paying a slight \GCH\ penalty, we can employ the diagonalization argument to obtain a completely local version:

\theorem Completely Local Version.(Woodin) If $\k$ is $\l$-supercompact and $2^{\l^\ltk}\!=\l\plus$, then this is preserved by fast function forcing. Indeed, every $\l$-super\-compact\-ness embedding in the ground model lifts to the forcing extension.
\ref\LocalB

\proof Suppose $j:V\to M$ is a $\l$-supercompact embedding in $V$, and that $f$ is a $V$-generic fast function. Let $p$ be the condition $\sing{\<\k,\l>}$, so that below $p$ the forcing $j(\F)$ factors as $\F*\Ftail$ and $\Ftail$ is $\ltel$-closed in $M$. Since a simple counting argument shows $\card{j(2^\k)}=\l\plus$, there are at most $\l\plus$ many open dense subsets of $\Ftail$ in $M$, counted in $V$.  Thus, using the closure of $M$ and the closure of the forcing, I may line them up and diagonalize against them to construct in $V$ an $M$-generic filter $\ftail\of\Ftail$. Consequently, in $V[f]$ the embedding lifts to $j:V[f]\to M[j(f)]$ where $j(f)=f\union p\union\ftail$, and it is not difficult to verify that this embedding is a $\l$-supercompact embedding in $V[f]$.\qed

Because the previous theorems show that fast function forcing preserves large cardinals, one expects many embeddings $j:V[f]\to M[j(f)]$ in the fast function extension. What is more---and this is the fundamental fact which makes fast function forcing useful---the next theorem shows that these embeddings are so easily modified that the value of $j(f)(\k)$ can be almost arbitrarily specified. Let me define that a measure $\eta$ in $V[f]$, or any forcing extension, is {\df standard} when the critical point $\k$ of the induced embedding $j:V[f]\to M[j(f)]$ is definable in $M[j(f)]$ from $s=[\id]_\mu$ and parameters in $\ran(j\restrict V)$. Thus, any normal measure on $\k$ is standard, as is any supercompactness measure (since $\k$ is the least element not in $j\image\l$). Also, Lemma 2.7 below shows that in the type of forcing extensions of this paper, every $\theta$-strong compactness measure is isomorphic to a standard strong compactness measure.

\theorem Fast Function Flexibility Theorem. Suppose that $f\from\k\to\k$ is a fast function added generically over $V$ and that $j:V[f]\to M[j(f)]$ is an embedding (either internal or external to $V[f]$) with critical point $\k$. Then for any $\a<j(\k)$ there is another embedding $j^*:V[f]\to M[j(f)]$ such that:
\points 1. $j^*(f)(\k)=\a$,\cr
            2. $j^*\restrict V=j\restrict V$,\cr
            3. $M[j^*(f)]\of M[j(f)]$, and\cr
            4. If $\a$ is not too much larger than $\k$ (see below), then $M[j^*(f)]=M[j(f)]$. In this case, if $j$ is the ultrapower by a standard measure $\eta$ concentrating on a set in $V$, then $j^*$ is the ultrapower by a standard measure $\eta^*$ concentrating on the same set, and moreover $\eta\intersect V=\eta^*\intersect V$ and $[\id]_\eta=[\id]_{\eta^*}$.\cr

\ref\Flexibility

\proof Fix $j$ and $\a$. Let $\bar\k$ be the next inaccessible above both $\k$ and $\a$, and let $\g$ be the next element of $\dom(j(f))$ above $\bar\k$. Thus, $\g$ is not a limit of inaccessible cardinals. By the Fast Function Factor Lemma, $\ftail=j(f)\restrict[\g,j(\k))$ is $M$-generic for $\F_{\g,j(\k)}$. 
Now consider the embedding $j\restrict V:V\to M$ (which perhaps may not be definable in $V$), and the condition $p=\sing{\<\k,\a>,\<\bar\k,\b>}$ where $\b<\g$ is larger than every inaccessible below $\g$. $$\vbox to1.5in{}$$ Below this condition, the forcing $j(\F)$ factors as $\F\cross\F_{\g,j(\k)}$. Since we have $M$-generics for these posets, we can let $j^*(f)=f\union p\union\ftail$ and lift the embedding to $j^*:V[f]\to M[j^*(f)]$. By construction we have $j^*(f)(\k)=p(\k)=\a$ and $j^*\restrict V=j\restrict V$. Also, $j^*(f)$ is easily constructed from $j(f)$, so $M[j^*(f)]\of M[j(f)]$. Finally, in the case that $\a$ does not exceed the next inaccessible cluster point of $\dom(j(f))$ beyond $\k$, then it follows that the `missing' part of $j(f)$, namely $j(f)\restrict[\k,\g)$, is simply a condition in $j(\P)$, and there lies in $M$. In this case $M[j(f)]=M[j^*(f)]$.

Finally, suppose in this case that $j$ is the ultrapower by a standard measure $\eta$ concentrating on a set $D\in V$. Let $s=[\id]_\eta$. This is a seed for $\eta$ in the sense that $X\in\eta\iff s\in j(X)$. Since $D\in\eta$ it follows that $s\in j(D)$ and consequently $s\in M$. For a technical reason, I will choose $\b$ in the previous argument to be an index of the condition $j(f)\restrict[\k,\g)\in M$ with respect to $j(\vec a)$ where $\vec a$ is an enumeration of $V_\k$ in $V$ such that for every $\xi<\k$ every element of $V_\xi$ appears unboundedly often among the first $\beth_\xi$ many elements of $\vec a$. Let $\eta^*$ be the measure germinated by the seed $s$ via $j^*$; i.e. $X\in\eta^*\iff s\in j^*(X)$. In order to argue that $j^*$ is the ultrapower by $\eta^*$, it suffices to show that every element of $M[j^*(f)]$ has the form $j^*(h)(s)$ for some $h\in V[f]$ (see \cite[Ham98] for an elementary introduction to these seed techniques). Let $\X$ be the seed hull of $s$, that is, the set of the elements in $M[j^*(f)]$ having this form. It is easy to verify the Tarski-Vaught criterion, and so $\X\elesub M[j^*(f)]$. Furthermore, since the measure $\eta$ was standard, it follows that $\k\in\X$ (and this is the only reason for that assumption). Consequently, $\b\in\X$ and so by the technical choice of $\b$ the missing part of $j(f)$ also lies in $\X$. Thus, from $j^*(f)\in\X$ we can reconstruct $j(f)$, and so $j(f)\in\X$. Now, suppose $x\in M[j^*(f)]=M[j(f)]$. Since $j$ is the ultrapower by $\eta$ we know that $x=j(h)(s)$ for some function $h\in V[f]$. This function has a name $\hdot\in V$. So $x=j(\hdot_f)(s)=j(\hdot)_{j(f)}(s)$. Since all the sets in this last expression are in $\X$, it must be that $x\in\X$ also; so $j^*$ is the ultrapower by $\eta^*$. The rest of the theorem follows because $[\id]_{\eta^*}=s=[\id]_\eta$ and $j^*\restrict V=j\restrict V$.\qed

In the context of a strongly compact cardinal $\k$, Menas was very concerned in \cite[Men74] with the situation in which there is a function $f\from\k\to\k$ with what I will call the {\df Menas property}, namely, that for every $\l$ there should be a fine measure $\mu$ on $P_\k\l$ with ultrapower embedding $j:V\to M$ such that $j(f)(\k)\geq\card{[\id]_\mu}^{M}$. These functions figured crucially in his preservation arguments.  Menas proved that every strongly compact limit of strongly compact cardinals has such a function, but conjectured that this would not be the case for every strongly compact cardinal. I will prove here, however, that one can have such a function for any strongly compact cardinal.

\theorem. Every fast function on a strongly compact cardinal has the Menas property.
\ref\STRCMenas

\proof This almost follows directly from the Flexibility Theorem, except for the difficulty that for large $\a$ the embedding $j^*$ produced in the Flexibility Theorem may not itself be a $\l$-strong compactness embedding; so an additional factor argument is needed. Suppose that $j:V[f]\to M[j(f)]$ is a $\l$-strongly compact embedding by some measure $\mu$ in $V[f]$. Let $\g=\card{s}$ where $s=[\id]_\mu$. Since $s\in M[j(f)]$ and $s$ has size $\g$ it follows that $s\in M[j(f)\restrict \g]$, and so it has a name $\dot s\in M$ of size $\g$. Since $j\image\l\of s$ by fineness, we may use the name $\dot s$ to build a set $\tilde s\in M$ of size $\g$ such that $j\image\l\of\tilde s$. Furthermore, we may assume $\k$ is the least element not in $\tilde s$, by simply removing it if necessary. Now let $j^*:V[f]\to M[j^*(f)]$ be an embedding as in the Flexibility Theorem such that $j^*(f)(\k)=\a$ for some $\a>\g$. Let $\tilde\mu$ be the measure germinated by the seed $\tilde s$ via $j^*$, so that $X\in\tilde\mu\iff\tilde s\in j^*(X)$. Since $j^*\image\l=j\image\l\of\tilde s\in j^*(P_\k\l)$, this measure is a fine measure on $P_\k\l$, and since it was obtained by a seed via $j^*$, we obtain the following factor diagram: $$\trianglediagram{V[f]}{j_0}{j^*}{\tilde M[j_0(f)]}{k}{M[j^*(f)]}$$ where $j_0$ is the ultrapower by $\tilde\mu$ and $k$ is the inverse collapse of the seed hull $\X=\set{j^*(h)(\tilde s)\st h\in V[f]}\elesub M[j^*(f)]$. Since $\k$ is the least element not in $\tilde s$, it follows that $\k\in\X$ and hence also $\a=j^*(f)(\k)\in\X$. Let $s_0$ and $\a_0$ be the collapses of $\tilde s$ and $\a$, respectively, so that $k(s_0)=\tilde s$ and $k(\a_0)=\a$. It follows that $[\id]_{\tilde\mu}=s_0$ and $j_0(f)(\k)=\a_0>\card{s_0}$, so $f$ has the Menas property with respect to $\tilde\mu$, as desired.\qed

The Menas property has a natural analogue for supercompact and strong cardinals. Specifically, I define for a supercompact cardinal $\k$ that $f\from\k\to\k$ has the {\df supercompact Menas property} when for every $\l$ there is a $\l$-supercompactness embedding $j$ for which $j(f)(\k)>\l$. Thus, for example, every Laver function has the Menas property. For a strong cardinal $\k$, I define that $f$ has the {\df strong Menas property} when for every $\l$ there is a $\l$-strong embedding $j$ for which $j(f)(\k)>\beth_\l$. Such functions are related to the high-jumping functions of \cite[Ham98].

\theorem. Every fast function on a supercompact cardinal has the supercompact Menas property. 
\ref\SCMenas

\proof This theorem is true level-by-level for partially supercompact cardinals. Suppose that $j:V[f]\to M[j(f)]$ is a $\l$-supercompactness embedding in $V[f]$, the ultrapower by a normal fine measure $\eta$ on $P_\k\l$. By Remark \FFFGap, we know that $M\of V$ and in fact $j:V\to M$ is definable in $V$ and $M$ is $\l$-closed there (though it need not be the ultrapower by a normal measure on $P_\k\l$ there). I claim that $j(f)\restrict[\k,\l)$ is in $M$. If not, then part of it must be generic over $M$ for some nontrivial $\lte\k$-closed forcing of size at most $\l$, namely, $\F_{\k,\g}$, where $\g$ is the first inaccessible cluster point of $\dom(j(f))$ beyond $\k$. Since this poset is the same in $M$ as in $V$, with the same dense sets, the filter generated by $j(f)\restrict[\k,\g)$ in $\F_{\k,\g}$ must be $V$-generic. But this is impossible, since it was added by the $\k$-c.c. forcing $\F$. So $j(f)\restrict[\k,\l)\in M$. Consequently, by the Flexibility Theorem, we may modify the embedding to $j^*:V[f]\to M[j^*(f)]$ so that $j^*(f)(\k)=\a$ for some $\a>\l$ and $M[j^*(f)]=M[j(f)]$. Furthermore, we may assume that $j^*$ is the ultrapower by a measure $\eta^*$ with $[\id]_\eta=[\id]_{\eta^*}$. Since $\eta$ is normal and fine, it follows that $[\id]_\eta=j\image\l$, so $[\id]_{\eta^*}=j\image\l=j^*\image\l$, and so $\eta^*$ is also a normal fine measure on $P_\k\l$. Finally, since $j^*(f)(\k)=\a>\l$, the measure $\eta^*$ exhibits that $f$ has the Menas property for a $\l$-supercompactness embedding.\qed

The previous proof in fact shows that the function $h:\k\to\k$, where $h(\g)$ is the next inaccessible cluster point of $\dom(f)$ beyond $\g$, is a high-jumping function in the terminologoy of \cite[Ham98]. It follows, by Theorem 3.4 of \cite[Ham98], that fast function forcing must destroy the almost hugeness of $\k$. 

\theorem. Every fast function on a strong cardinal has the strong Menas property. 
\ref\STRMenas

\proof This theorem is almost true level-by-level. Specifically, I will show that if $\k$ is $(\l+1)$-strong in $V[f]$ then $f$ has the Menas property with respect to $(\l+1)$-strong embeddings in $V[f]$. It follows, using the usual factor argument and the induced $\l$-strong extender, that $f$ has the Menas property with respect to a $\l$-strong embedding also. So, suppose $j:V[f]\to M[j(f)]$ is a $(\l+1)$-strong embedding, so that $V_{\l+1}[f]\of M[j(f)]$. By factoring through by the natural extender I may assume that $M[j(f)]$ is closed under $\k$-sequences in $V[f]$, and consequently, by Remark \FFFGap, that $M\of V$ and furthermore $M_{\l+1}=V_{\l+1}$. I will argue as in the previous theorem that $j(f)\restrict[\k,\l)\in M$. If this fails, then there must be some $\g\leq\l$ which is an inaccessible cluster point of $\dom(j(f))$, and $j(f)\restrict[\k,\g)$ is $M$-generic for for $\F_{\k,\g}^M$. Since $\F_{\k,\g}$ is the same whether computed in $V$ or $M$, and $V$ and $M$ have the same dense sets for it, it follows that $j(f)\restrict[\k,\g)$ is actually $V$-generic for $\F_{\k,\g}$. But this is impossible since $j(f)\in V[f]$, a $\k\plus$-c.c. forcing extension of $V$, and $\F_{\k,\g}$ is $\ltek$-closed. So $j(f)\restrict[\k,\g)$ must just be a condition in $j(\F)$ and hence an element of $M$. Now, we continue as in the previous theorem. By the Flexibility theorem, there is another embedding $j^*:V[f]\to M[j^*(f)]$ such that $j^*(f)(\k)>\beth_{\l+1}$ and $M[j^*(f)]=M[j(f)]$. There is no trouble making $j^*(f)(\k)$ larger than $\beth_{\l+1}$ since the proof of the Flexibility Theorem shows that it can easily be pushed up beyond the next inaccessible above $\l$. Thus, $j^*:V[f]\to M[j^*(f)]$ has $(V[f])_{\l+1}\of M[j^*(f)]$ and $j^*(f)(\k)>\beth_{\l+1}$, as desired.\qed

Let me conclude this section with a quick application of fast function forcing. Kunen and Paris \cite[KunPar71] were the first to show that a measurable cardinal $\k$ can have many normal measures in a forcing extension. The following argument shows that fast function forcing works nicely to see the same fact for a variety of large cardinals. 

\theorem Many Measures Theorem. Fast function forcing with $\k$ adds many measures. Specifically, 
\points	1. Every (sufficiently nice) weak compactness filter on $\k$ in $V$ extends to $\k$ many weak compactness filters in $V[f]$.\cr
	2. If $2^\k=\k\plus$, then every measure on $\k$ in $V$ extends to $2^{2^\k}$ many measures in $V[f]$, the maximum conceivable number. Indeed, every measure in $V$ extends to $2^{2^\k}$ many measures, each isomorphic in $V[f]$ to a distinct normal measure.\cr
	3. If $\k$ is $2^{\l^\ltk}\!$-strongly compact in $V$ then there are $\l\plus$ many non-isomorphic $\l$-strong compactness measures in $V[f]$. Thus, if also $2^\k=\k\plus$, then there are $2^{2^\k}\!\cdot\l\plus$ many.\cr
	4. If $2^{\l^\ltk}\!=\l\plus$ then every $\l$-supercompactness measure in $V$ extends to $2^{2^{\l^\ltk}}$ many $\l$-supercompactness measures in $V[f]$, the maximum conceivable number.\cr

\proof Though it is a bit more work to get the optimal bounds, this theorem follows in spirit from the Flexibility Theorem; essentially, the fact that $j(f)(\k)$ can have many different values means that there must be many different measures. Thus, 1 holds for the nice filters I managed to lift in the previous theorems, because for each weak compactness embedding there are $\k$ many possible values for $j(f)(\k)$.

Let me prove 2. The simple idea of looking at the possible values of $j(f)(\k)$ easily gives $2^\k$ many measures; in order to get $2^{2^\k}$ many measures, I will consider the possible values of $j(f)$. Suppose $\k$ is measurable in $V$ and $\mu$ is any measure on $\k$ with embedding $j_\mu:V\to M$. Fast Function Theorem \FFMeasurable\ shows that there are many lifts of $j_\mu$ to $j:V[f]\to M[j(f)]$. How many are there? Well, the proof proceeded by diagonalizing against the dense sets of $M$, and there are diverse ways to carry out this diagonalization. Specifically, below any condition in $j(\F)$ there is an antichain of size $j(\k)$, which has size $\k\plus$ in $V$. Thus we can build a tree of height $\k\plus$ of descending conditons in $\Ftail$ such that every node splits into an antichain of size $\k\plus$ on the next level. Furthermore, we can arrange that every node on the $\a^\th$ level of this tree is in the $\a^\th$ dense set of $M$, so that any $\k\plus$-branch through this tree will produce an $M$-generic for $\Ftail$. Since there are ${\k\plus}^{\k\plus}=2^{2^\k}$ many $\k\plus$-branches through this tree, there are $2^{2^\k}$ many ways to perform the diagonalization, and each of the resulting generics produces a different $j(f)$, and consequently a different measure in $V[f]$. So we have many measures in $V[f]$. Now, let me argue that we can arrange for all of these embeddings to be ultrapowers by normal measures in $V[f]$. If we build the tree below the condition which ensures $j(f)(\k)=\a$, where $\a=[\id]_\mu$ is the canonical seed for $\mu$, then with respect to the lifted embedding $j:V[f]\to M[j(f)]$, the seed hull of $\k$ generates the old seed $\a$ and consequently {\it all} of $M[j(f)]$ (see the Old Seed Lemma of \cite[Ham97]). Thus, since the entire embedding is in the seed hull  of $\k$, the embedding is an embedding by the normal measure $\eta$ induced by $\k$. If $\nu$ is the measure germinated by the seed $\a$ with respect to $j$, then $\mu$ extends to $\nu$ since we lifted the embedding. But since $\k$ generates $\a$ and vice versa, the measures $\nu$ and $\eta$ are isomorphic, so 2 holds.

Statement 4 holds similarly. Suppose that $2^{\l^\ltk}\!=\l\plus$ and $j:V\to M$ is a $\l$-supercompactness embedding. The diagonalization technique of \LocalB\ shows that we may lift the embedding to $j:V[f]\to M[j(f)]$. Again, we can build a $\ltel$-closed tree of height $\l\plus$ and $\l\plus$ branching at each node such that any branch through this tree provides a different generic $j(f)$ with which to lift the embedding. By the Old Seed Lemma of \cite[Ham97], the seed $j\image\l$ still generates the whole embedding, and consequently each of these lifts provides a different $\l$-supercompactness measure lifting and extending the original measure. Thus, there are $(\l\plus)^{\l\plus}=2^{2^{\l^\ltk}}$ many measures extending the original measure, as desired. By working below a condition which forces $j(f)(\k)=\l+1$, we can arrange that all these measures witness the Menas property of $f$. 

Finally, let me prove 3. By \STRCMenas, there are fine measures $\mu$ on $P_\k\l$ in $V[f]$ witnessing the Menas property of $f$, so that the corresponding embedding $j:V[f]\to M[j(f)]$ has $j(f)(\k)>\card{s}$ where $s=[\id]_\mu$. Furthermore, we may assume that $\k$ is the least element not in $s$, by simply removing it if necessary and working with the induced isomorphic measure; so we may assume that the measure $\mu$ is standard. Thus, by statement 4 in the Flexibility Theorem, we may for any $\a<\l\plus$ find an embedding $j^*:V[f]\to M[j^*(f)]$, the ultrapower by a measure $\mu^*$ with $[\id]_{\mu^*}=s$, such that $j^*(f)(\k)=\a$, $j^*\restrict V=j\restrict V$ and $M[j^*(f)]=M[j(f)]$. Since $j^*\image\l=j\image\l\of s\in j^*(P_\k\l)$, it follows that $\mu^*$ is a fine measure on $P_\k\l$. And since different choices of $\a$ provide different embeddings $j^*$, these measures are all pairwise non-isomorphic, and so we have $\l\plus$ many measures. Finally, since the argument before Theorem 4.2 shows that the product of a normal measure with a strong compactness measure is isomorphic to a strong compactness measure, by 2 if $2^\k=\k\plus$ there are at least $2^{2^\k}$ many strong compactness measures on $P_\k\l$ in $V[f]$, and so the theorem is proved.\qed

Previously, it was not known even how to force two non-isomorphic $\l$-strong compactness measures for a strongly compact cardinal. Nevertheless, in the case of strong compactness, the theorem is not the strongest conceivable result, since the following question remains open.

\question. Suppose $\k$ is strongly compact. Is there a forcing extension in which for every $\l$ there are the maximum conceivable number of non-isomorphic fine measures on $P_\k\l$, namely $2^{2^{\l^{{<}\k}}}\!\!$ many?
 
\section Generalized Laver Functions

The existence of Laver functions in the supercompact cardinal context has proved indispensible; these functions appear in dozens if not hundreds of articles. Because of this, we would really like to have Laver functions for other kinds of large cardinals. I am pleased, therefore, to prove here that fast function forcing adds a new completely general kind of Laver function to any large cardinal, thereby freeing the notion of Laver function from the supercompact cardinal context.

Specifically, I define in $\Vbar$ that $\ell\from\k\to\Vbar_\k$ is a {\df generalized Laver function} under the function $f\from\k\to\k$ when for any embedding $j:\Vbar\to\Mbar$ with critical point $\k$ and any $z\in H(\l\plus)^\Mbar$ where $\l=j(f)(\k)$ there is another embedding $j^*:\Vbar\to\Mbar$ (yes, the same $\Mbar$) such that $j^*(\ell)(\k)=z$ and $j^*\restrict\ORD=j\restrict\ORD$. This definition is perfectly sensible whether $\k$ is measurable, strong, strongly compact, supercompact, or huge, and so on. Of course, I make this definition only in the nontrivial case that $j(f)(\k)>0$ is possible; naturally the function $\ell$ would be the most useful when the function $f$, like a fast function, has the property that $\l=j(f)(\k)$ can be very large. Often, the function $f$ will in fact be a fast function or at least have the Menas property. In this case, every generalized Laver function on a supercompact cardinal is a Laver function in Laver's original sense. But the converse need not hold; indeed, there may be no generalized Laver functions at all:

\theorem Observation. If $V=\HOD$ and $\k$ is at least measurable, then there is no generalized Laver function for $\k$. In particular, there is no generalized Laver function in $L[\mu]$ or in the core models.

\proof Suppose $V=\HOD$ and $\k$ is measurable. Let $j:V\to M$ be any embedding with critical point $\k$. If $j^*:V\to M$ is an embedding such that $j\restrict \ORD=j^*\restrict \ORD$ then since every set is hereditarily ordinal definable, it follows that $j=j^*$.  Consequently, there is no freedom to choose $j^*(\ell)(\k)$; it must be equal to
$j(\ell)(\k)$. So there is no generalized Laver function.\qed Observation

\theorem Generalized Laver Function Theorem. Fast function forcing adds a generalized Laver function. Specifically, after fast function forcing $V[f]$, there is a function $\ell\from\k\to (V[f])_\k$ with the property that for any embedding $j:V[f]\to M[j(f)]$ with critical point $\k$ (whether internal or external) and for any $z\in H(\l\plus)^{M[j(f)]}$, where $\l=j(f)(\k)$, there is another embedding $j^*:V[f]\to M[j(f)]$ such that: 
\points 1. $j^*(\ell)(\k)=z$, \cr
        2. $M[j^*(f)]=M[j(f)]$,\cr
        3. $j^*\restrict V=j\restrict V$, and\cr
        4. If $j$ is the ultrapower by a standard measure $\eta$ concentrating on a set in $V$, then $j^*$ is the ultrapower by a standard measure $\eta^*$ concentrating on the same set and moreover $\eta^*\intersect V=\eta\intersect V$ and  $[\id]_{\eta^*}=[\id]_\eta$.\cr

\ref\LaverFunction 

\proof The idea is quite simple, given the Flexibility Theorem for fast function forcing. In $V$ enumerate $V_\k$ as $\vec a=\<a_\a\st\a<\k>$ with the property that for every $\xi<\k$ every element of $V_\xi$ appears unboundedly often among the first $\beth_\xi$ many elements of the enumeration. In $V[f]$ let $\ell(\g)=(a_{f(\g)})_{f\restrict\g}$, provided that this makes sense, i.e., that $\g\in\dom(f)$ and $a_{f(\g)}$ is an $\F_\g$-name. Suppose $j:V[f]\to M[j(f)]$ is given and $z\in H(\l\plus)^{M[j(f)]}$ where $\l=j(f)(\k)$. By the closure of the tail forcing, $z\in M[f]$ and so $z=\dot z_f$ for some name $\dot z\in M$. The name $\dot z$ must be $j(\vec a)(\a)$ for some index $\a$, and by the assumption on $\vec a$ such an $\a$ can be found below the next inaccessible beyond $\l$ and $\k$. Therefore, by the Flexibility Theorem, there is another embedding $j^*:V[f]\to M[j^*(f)]=M[j(f)]$ satisfying the conclusions of the Flexibility Theorem, with $j^*(f)(\k)=\a$. In particular, statements 2, 3 and 4 hold. It follows, by the definition of $\ell$, that $j^*(\ell)(\k)=\dot z_f=z$, as desired for statement 1. So the theorem is proved.\qed

Notice that the restriction that $z$ is in $H(\l\plus)$ is not onerous, because by the Flexibility Theorem the value of $\l=j(f)(\k)$ is highly mutable and can be made to be any desired ordinal up to $j(\k)$. Certainly, any $z$ in $H(\d)^{M[j(f)]}$ can be accomodated without modification for any $\d$ up to the next inaccessible cluster point of $\dom(j(f))$ beyond $\k$. And the arguments of \SCMenas\ and \STRMenas\ show that for $\l$-supercompact or $(\l+1)$-strong embeddings, this is always at least $\l$. More generally, though, any element of $M_{j(\k)}[f]$ is a possible value of $j^*(\ell)(\k)$, because given any $j:V[f]\to M[j(f)]$ one can first apply the Flexibility Theorem to get $j^*:V[f]\to M[j^*(f)]\of M[j(f)]$ such that $j^*(f)(\k)$ is large, and then apply the Generalized Laver Function Theorem to make $j^*(\ell)(\k)$ whatever element of $M_{j(\k)}[f]$ was desired.

For the remainder of this section, let me simply spell out the particular consequences of the previous theorem for various large cardinals. Henceforth in this section, therefore, let $\ell$ be the generalized Laver function of the previous theorem computed in the fast function extension $V[f]$ relative to the fixed enumeration $\vec a$ of $V_\k$. 

\theorem. Suppose in $V$ that $\k$ is measurable and $2^\k=\k\plus$. Then for any ultrapower embedding $j:V\to M$ by a measure on $\k$ and any $z\in H(\k\plus)^{V[f]}$ there is a lift $j:V[f]\to M[j(f)]$ such that $j(\ell)(\k)=z$. Furthermore, the lift can be arranged to be a normal ultrapower.

\proof This is what falls out of the previous arguments. Beginning with any ultrapower $j:V\to M$ by the measure $\mu$ and any $z\in H(\k\plus)^{V[f]}$ it follows that $z\in M[f]$ and so we can lift the embedding to $j:V[f]\to M[j(f)]$ in such a way that $j(f)(\k)$ picks out the ordinal index of a name for $z$, so that $j(\ell)(\k)=z$. By using a name, say, which also codes the ordinal $\a=[\id]_\mu$, it follows that for some function $g\in V[f]$ we have $\a=j(g)(\k)$. Since every element of $M$ has the form $j(g')(\a)$ where $g'$ is a function in $V$, it follows from this that every element of $M[j(f)]$ has the form $j(h)(\k)$ for some function $h\in V[f]$, and consequently the lifted embedding is a normal ultrapower.\qed

Essentially the same argument works for weakly compact cardinals:

\theorem. Suppose in $V$ that $\k$ is weakly compact. Then for any set $z\in H(\k\plus)$ in $V[f]$ there is a weakly compact embedding $j:N[f]\to M[j(f)]$ such that $j(\ell)(\k)=z$. In particular, $\Diamond_\k$ holds. 

\proof Since $z\in V[f]$ there is a name $\dot z\in V$ of hereditary size $\k$. Pick a transitive $N\elesub H(\k\plus)$ of size $\k$ in $V$ such that $\dot z, \F,\dot\ell,\vec a\in N$ and $N$ is closed under $\ltk$-sequences. Since $\k$ is weakly compact in $V$ there is an embedding $j:N\to M$ with critical point $\k$. As in \FFFWC, we can assume that $M$ is also closed under $\ltk$-sequences. The usual argument shows $P(\k)^N\of M$ and so $\dot z\in M$. By the diagonalization argument, because there are only $\k$ many dense subsets of $\F$ in $N$, we can lift the embedding to $j:N[f]\to M[j(f)]$, and furthermore, we can do so in such a way that $j(f)(\k)=\a$ where $\a$ is the index of $\dot z$ with respect to $j(\vec a)$. Consequently, $j(\ell)(\k)=\dot z_f=z$, as desired. It is easy now to deduce that $f$ is a powerful kind of $\Diamond$ sequence.\qed

Let me now gradually move upwards through the large cardinal hierarchy.

\theorem. If $j:V[f]\to M[j(f)]$ is a $\l$-strong embedding (with the natural extender) and $z\in (V[f])_\l$ then there is another $\l$-strong embedding $j^*:V[f]\to M[j^*(f)]$ such that $j^*(\ell)(\k)=z$, $j^*\restrict V=j\restrict V$ and $M[j^*(f)]=M[j(f)]$. 

\proof Suppose that $j:V[f]\to M[j(f)]$ is a $\l$-strong embedding generated by the natural extender, so that $M[j(f)]=\set{j(h)(s)\st h\in V[f]\and s\in (V[f])_\l}$, and that $z\in (V[f])_\l$. If $\l$ is a successor ordinal then the argument of \STRCMenas\ shows that $j(f)\restrict[\k,\l)\in M$, and so $z=\dot z_f$ for some $\dot z\in M$. Alternatively, if $\l$ is a limit ordinal, then $z\in (V[f])_\b$ for some much smaller $\b$, and consequently the argument of \STRMenas\ applied to the induced factor embedding shows $j(f)\restrict[\k,\b)\in M$. Thus, $z=\dot z_f$ for some $\dot z\in M_\l$. In either case, the name $\dot z$ has, below the next inaccessible, some index $\a$ with respect to $j(\vec a)$, and so by the Flexibility Theorem, there is another embedding $j^*:V[f]\to M[j^*(f)]$ such that $j^*(f)(\k)=\a$, $j^*\restrict V=j\restrict V$ and $M[j^*(f)]=M[j(f)]$. By the choice of $\a$ it follows that $j^*(\ell)(\k)=z$, and so the theorem is proved.\qed

As before, the theorem can be modified to allow for $z$ which appear higher in the hierarchy if we are willing to give up the equality of $M[j^*(f)]$ and $M[j(f)]$. Specifically, if $z\in M_{j(\k)}[f]$, then there will be an embedding $j^*:V[f]\to M[j(f)]$ such that $j^*(\ell)(\k)=z$, $j^*\restrict V=j\restrict V$ and $M[j^*(f)\of M[j(f)]$. The embedding $j^*$ will still be a $\l$-strong embedding because $(M[j(f)])_\l=(M[f])_\l$ by the argument showing $j(f)\restrict[\k,\b)\in M$. The next inaccessible cluster point of $\dom(j(f))$ beyond $\k$ must be at least $\l$. 

Next, I treat the case of strongly compact cardinals. 

\theorem. If the embedding $j$ of \LaverFunction\ is a $\theta$-strong compactness embedding, where $\theta^\ltk=\theta$, then the embedding $j^*$ may also be chosen to be a $\theta$-strong compactness embedding.
\ref\STRCLaverFunction

Since every $\theta$-strongly compact measure is in fact isomorphic to a $\theta^\ltk$-strongly compact measure, we see by simply replacing $\theta$ with $\theta^\ltk$ that the assumption that $\theta^\ltk=\theta$ is hardly a restriction at all. And because a measure $\mu$ is a $\theta$-strong compactness measure exactly when $s=[\id]_\mu$ is a cover of $j_\mu\image\theta$ with a subset of $j_\mu(\theta)$ of size less than $j_\mu(\k)$, the theorem follows by statement 4 of \LaverFunction\ and the following lemma. Define that a forcing extension $V[G]$ is {\df mild} if every set of hereditary size less than $\k$ in $V[G]$ is added by a poset of size less than $\k$ in $V$. Certainly fast function forcing is mild, because all the tail forcings $\F_{\l,\k}$ are $\ltel$-closed. 

\lemma. Every $\theta$-strong compactness measure in a mild forcing extension $V[G]$, where $\theta^\ltk=\theta$, is isomorphic to a standard $\theta$-strong compactness measure which concentrates on $(P_\k\theta)^V$.
\ref\STRCMild

\proof Suppose that $j:V[G]\to M[j(G)]$ is the ultrapower by a fine measure $\eta$ on $P_\k\theta$ in $V[G]$. 
Since measures are isomorphic exactly when they induce the same embedding (see \cite[Ham97]), it suffices to show that $j$ is the ultrapower by a standard fine measure concentrating on $(P_\k\theta)^V$. Let $s=[\id]_\eta$. Thus, $j\image\theta\of s$ and $\card{s}^{M[j(G)]}<j(\k)$. Thus, by mildness, $s\in M[\Gtilde]$ for some generic $\Gtilde\of\Ptilde$ for forcing of some size $\g$  such that $\card{s}\leq\g<j(\k)$. Thus, $s$ has a name $\dot s\in M$ of size $\g$. Using this name it is possible to construct a set $\tilde s\in M$ of size $\g$ such that $j\image\theta\of\tilde s$ and $\k\notin\tilde s$. Since $\theta^\ltk=\theta$, the measure $\eta$ is isomorphic to a $\k$-complete measure on $\theta$. There must therefore be an ordinal $\d<j(\theta)$ such that $M[j(G)]=\set{j(h)(\d)\st h\in V[G]}$. We may assume, by simply adding such a point if necessary, that the largest element of $\tilde s$ has the form $\<\b,\d>$, using a suitable definable pairing function, for some ordinal $\b<j(\theta)$. Let $\tilde\eta$ be the measure germinated by $\tilde s$ via $j$. Since $\tilde s$ is a subset of $j(\theta)$ of size $\g<j(\k)$ and $j\image\theta\of\tilde s$, it follows that $\tilde\eta$ is a fine measure on $P_\k\theta$ in $V[G]$. Furthermore, since $\tilde s\in M$, it concentrates on the $P_\k\theta$ of the ground model $V$. I claim that $\eta$ is isomorphic to $\tilde\eta$. To prove this, it suffices by the seed theory of \cite[Ham97] to show that the seed hull of $\tilde s$, namely $\X=\set{j(h)(\tilde s)\st h\in V[G]}\elesub M[j(G)]$, is all of $M[j(G)]$. By the assumption on the largest element of $\tilde s$, we know $\d\in\X$ and since also $\ran(j)\of\X$, it follows that $M[j(G)]\of\X$, as desired. The measure $\tilde\eta$ is standard because $\k$ is the least element not in $\tilde s=[\id]_{\tilde\eta}$. So I have proved that every $\theta$-strong compactness measure $\eta$ in $V[G]$ is isomorphic to a standard $\theta$-strong compactness measure $\tilde\eta$ in $V[G]$ concentrating on the $P_\k\theta$ of the ground model. \qed

So Theorem \STRCLaverFunction\  is proved. As usual, by giving up the equality of $M[j^*(f)]$ and $M[j(f)]$ it is possible to accomodate larger $z$ than stated in the theorem, as I will prove next. Define that an embedding $j:\Vbar \to \Mbar $ has the $\theta$-strong compactness {\df cover property} when there is a set $s\in\Mbar $ such that $j\image\theta\of s$ and $\card{s}^{\Mbar }<j(\k)$. Thus, $s$ can be used to germinate via $j$ a fine measure on $P_k\theta$. In the event that $j$ is an ultrapower by a measure on some set, it follows by an easy argument that every subset of $\Mbar $ of size at most $\theta$ is covered by an element of $\Mbar$ of size $\card{s}^{\Mbar }$. 

\theorem. If $j:V[f]\to M[j(f)]$ is a $\theta$-strong compactness embedding and $z\in M_{j(\k)}[f]$, then there is another embedding $j^*:V[f]\to M[j^*(f)]$ with the $\theta$-strong compactness cover property such that $j^*(\ell)(\k)=z$,  $M[j^*(f)]\of M[j(f)]$ and $j^*\restrict V=j\restrict V$. 

\proof Suppose $j:V[f]\to M[j(f)]$ is a $\theta$-strong compactness embedding in $V[f]$ and $z\in M_{j(\k)}[f]$. Thus, $z=\dot z_f$ for some name $\dot z\in M_{j(\k)}$. By the Flexibility Theorem, I may find $j^*:V[f]\to M[j^*(f)]$ such that $j^*\restrict V=j\restrict V$, $M[j^*(f)]\of M[j(f)]$ and, most importantly, $j^*(f)(\k)$ picks out the index of $\dot z$ with respect to $j(\vec a)$, so that $j^*(\ell)(\k)=(\dot z)_f=z$. By Lemma \STRCMild\ we may assume that there is a set $s\in M$ such that $j\image\theta\of s$ and $\card{s}^M<j(\k)$. This set also works as a cover, therefore, in $M[j^*(f)]$. \qed

\theorem. If the embedding $j$ of \LaverFunction\ is a $\theta$-supercompactness embedding, then the embedding $j^*$ may also be chosen to be a $\theta$-supercompactness embedding. 

\proof This is immediate by property 4 of \LaverFunction\ and Remark \FFFGap, since the remark shows that $j\image\theta\in M$, and this is $[\id]_\eta$. That is, supercompactness measures in $V[f]$ are always standard, and they always concentrate on the $P_\k\theta$ of the ground model $V$.\qed

The next theorem improves on this; even when $\l=j(f)(\k)$ is small  it is possible for $\theta$-supercompactness embeddings to have $j^*(\ell)(\k)=z$ for any $z\in H(\theta\plus)$.

\theorem. If $j:V[f]\to M[j(f)]$ is a $\theta$-supercompactness embedding in $V[f]$ and $z\in H(\theta\plus)^{V[f]}$ then there is another $\theta$-supercompactness embedding $j^*:V[f]\to M[j^*(f)]$ such that $j^*(\ell)(\k)=z$, $j^*\restrict V=j\restrict V$ and $M[j^*(f)]=M[j(f)]$. 

\proof This follows by the same idea as in \SCMenas; the point is that $j(f)\restrict[\k,\theta)$ must be a condition in $M$, and so the value of $j(f)(\k)$ can be freely changed so as to pick out the index of the name of any element in $H(\theta\plus)^{M[j(f)]}$.\qed

The previous theorem is true level-by-level in the sense that it is true even when $\k$ is only partially supercompact; for example, $\k$ may be only measurable. What's more, it is as before possible for the function $\ell$ to capture more $z$ than just those in $H(\theta\plus)$. Specifically, if $j:V[f]\to M[j(f)]$ is a $\theta$-supercompact embedding and $z\in M_{j(\k)}[f]$ then there is an embedding $j^*:V[f]\to M[j^*(f)]$ such that $j^*(\ell)(\k)=z$, $j^*\restrict V=j\restrict V$ and $M[j^*(f)]\of M[j(f)]$. In this case, however, the embedding $j^*$ may not be a $\theta$-supercompactness embedding, though by Remark \FFFGap\ it will have $j^*\image\theta=j\image\theta\in M$. 

\section The Lottery Preparation

I aim here to present the lottery preparation, a new general kind of Laver preparation, which works uniformly with a variety of large cardinals---such as weakly compact cardinals, measurable cardinals, strong cardinals, strongly compact cardinals and supercompact cardinals---and makes them indestructible by various further forcing, depending on the strength of the cardinal.

Let me begin by defining my terms. The basic building block is what I call a lottery sum. Specifically, the {\df lottery sum} of a collection $A$ of forcing notions is the forcing notion  $\oplus A=\set{\<\Q,p>\st \Q\in A\and p\in\Q}\union\sing{\one}$, ordered with $\one$ above everything and  $\<\Q,p>\leq\<\Q',q>$ when $\Q=\Q'$ and $p\leq_\Q q$. Because compatible conditions must have the same $\Q$, the forcing effectively holds a lottery among all the posets in $A$, a lottery in which the generic filter selects a `winning' poset $\Q$ and then forces with it.
$$\vbox{\vskip 1.5in}$$ Note that the lottery sum of the empty set is the trivial poset $\sing{\one}$. I will define the lottery preparation of $\k$ relative to a fixed function $f\from\k\to\k$. Though the definition works fine with any function, the forcing works best when used with a function having the Menas property, such as a fast function. 

The lottery preparation of $\k$ will be a $\k$-iteration which at many stages $\g<\k$ will perform the lottery sum of the collection of posets which are {\df allowed} at stage $\g$. Thus, at stage $\g$, the generic filter will effectively select a particular such poset as the winner of the lottery and then force with it. Generically, a wide variety of posets will be chosen in the lotteries below $\k$, thereby reflecting the possibilities at stage $\k$ on the $j$-side. The essential idea is that rather than consulting a Laver function about which particular forcing is to be done at stage $\g$, the lottery preparation instead uses the lottery sum of all posets which we might like to see at stage $\g$, and lets the generic filter decide generically amongst them.

Officially, let me say that a poset $\Q$ is {\df allowed} at stage $\g$ when for every $\d<\g$ the poset $\Q$ is $\ltd$-strategically closed (that is, the second player has a strategy enabling her to play a descending $\d$-sequence from the poset, where the players alternately play elements descending through the poset, and the second player plays at limit stages). This requirement, while broadly inclusive, is enough to ensure that the tail forcing is distributive. 

Let me now give the definition. The {\df lottery preparation} of $\k$ relative to the function $f\from\k\to\k$ is the reverse Easton support\footnote{*}{By reverse Easton support, I mean that it is a forcing iteration in which direct limits are taken at all inaccessible limit stages, and inverse limits at all other stages.} $\k$-iteration which has nontrivial forcing at stage $\g$ only when $\g\in\dom(f)$ and $f\image\g\of\g$. At such stages, the forcing $\Q_\g$ is the lottery sum in $V^{\P_\g}$ of all posets in $H(f(\g)\plus)$ which are allowed at stage $\g$. Otherwise, the forcing at stage $\g$ is trivial.

While I have proved in the previous section that fast function forcing adds a generalized Laver function, please observe that I am not using this generalized Laver function to define the lottery preparation. Certainly one could use the generalized Laver functions to define a kind of generalized Laver preparation, and such a preparation would have many of the same features (by essentially the same arguments) that I identify here for the lottery preparation. But it seems conceptually simpler to me, and more to the point, to use lottery sums in order to allow the generic filter to decide which forcing is to be done at each stage. Doing so avoids the need to carefully configure the embedding so that $j(\ell)(\k)$ is as required; with a lottery, one simply works below the condition which opts for the desired forcing at stage $\k$. Indeed, this ability to select arbitrarily the winner of the stage $\k$ lottery is what allows us to get by without any Laver function. In a sense the lottery preparation shows that what was truly important about Laver's function was not that $j(\ell)(\k)$ could be arranged to be any desired set---since a lottery sum can generically pick out any desired set at stage $\k$---but rather that the Laver function could be arranged so that the next element of the domain of $j(\ell)$ beyond $\k$ was as large as you like. This is what supports the crucial tail forcing arguments; as long as one has a way of reaching up high (e.g. by the Menas property), one can use a lottery sum to allow the generic filter to select any desired set or poset, and the tail forcing will be sufficiently closed. Thus, the lottery preparation defined relative to a fast function works effectively with a wide variety of large cardinals. 

\lemma Lottery Factor Lemma. For any $\g<\k$ which is closed under $f$, the lottery preparation $\P_\k$ factors as $\P_\g*\P_{\g,\k}$ where $\P_\g$ is the lottery preparation defined using $f\restrict\g$ and $\P_{\g,\k}$ is the lottery preparation defined in $V^{\P_\g}$ using $f\restrict[\g,\k)$. 

\proof This follows by the usual iterated forcing factor arguments. The point is that the $P_{\g+\a}$-names appearing in the stage $\g+\a$ lottery can be iteratively transformed, by recursion on $\a$, into $\P_\g*\P_{\g,\g+\a}$-names. There is no problem with the supports because we took an inverse limit at all but the inaccessible stages.\qed

\lemma. If in the lottery preparation there is no nontrivial forcing until beyond stage $\g$, then the preparation is $\lteg$-strategically closed.

\proof Since the forcing at each stage $\a>\g$ is $\a$-allowed, it is $\lteg$-strategically closed, with a (name of a) strategy $\sigma_\a$. Given a partial play, a descending sequence of conditions  $\<p^\b\st \b<\g'>$ for some $\g'<\g$, where $p^\b=\<\dot p^\b_\a\st \g<\a<\k>$, one applies the strategies $\sigma_\a$ coordinate-wise to obtain $\sigma(\<p^\b\st\b<\g'>)=\<\dot q_\a\st\g<\a<\k>$, where $\dot q_\a$ is the name for the condition obtained by applying the strategy $\sigma_\a$ to $\<\dot p^\b_\a\st\b<\g'>$. Recursively, since each of the strategies $\sigma_\a$ can successfully negotiate all the limits up to $\g$, so does this strategy $\sigma$, and so the lemma is proved.\qed

The consequence of this lemma is that when one factors the lottery preparation as $\P_\g*\P_{\g,\k}$, then $\P_{\g,\k}$ is $\ltg$-strategically closed in $V^{\P_\g}$. The two lemmas together show that if we are trying to lift an embedding $j:V\to M$ and $\Q$ is some forcing which which is allowed in the stage $\k$ lottery of $j(\P)$, then by simply working below a condition $p$ which opts for $\Q$ in the stage $\k$ lottery we may factor the forcing as $j(\P)\restrict p\iso\P*\Q*\Ptail$, where $\Ptail$ has trivial stages until beyond $j(f)(\k)$. Through this simple lottery technique, we obtain the crucial factorization that one ordinarily needs a Laver function to obtain, and we have done so in a completely general large cardinal context, with no supercompactness assumptions. This is the idea which will support the indestructibility results of the next section.

Let me now prove that the lottery preparation preserves a variety of large cardinals, beginning at the bottom and moving upwards. In the theorems below, if no assumption is explicity made concerning $f$, then the theorem holds for the lottery preparation defined using any function $f$. To avoid the triviality of small forcing, let me assume that the domain of $f$ is unbounded in $\k$. 

\theorem Lottery Preparation Theorem. The lottery preparation of an inaccessible cardinal $\k$ preserves the inaccessibility of $\k$. 

\proof Suppose that $\k$ becomes singular in $V[G]$. Let $\g\in\dom(f)$ be a closure point of $f$ larger than $\d=\cof(\k)^{V[f]}$ and factor the forcing at stage $\g$  as $\P_\g*\P_{\g,\k}$. The forcing $\P_\g$ has size less than $\k$ and so cannot have collapsed the cofinality of $\k$. The rest of the forcing $\P_{\g,\k}$ is $\lted$-strategically closed in $V^{\P_\g}$, and so also cannot have collapsed the cofinality of $\k$ to $\d$, a contradiction. Thus, $\k$ must be regular in $V[G]$. Similarly, the forcing $\P_\g$, since it has size less then $\k$, cannot force $2^\d\geq\k$ for any $\d<\k$, and the tail forcing $\P_{\g,\k}$ adds no further subsets to $\d$ if $\d<\g$. So $\k$ remains inaccessible in $V[G]$.\qed

\theorem Lottery Preparation Theorem. The lottery preparation of a weakly compact cardinal $\k$ preserves the weak compactness of $\k$. 

\proof Since the previous theorem shows that $\k$ remains inaccessible, it suffices to show that $\k$ has the tree property in $V[G]$. So, suppose $T$ is a $\k$-tree in $V[G]$. Choose a name $\dot T\in V$ so that $\dot T_G=T$. In $V$ pick a transitive $N\elesub H(\k\plus)$ of size $\k$ so that $\k,\dot T,\P\in N$ and $N^{\ltk}\of N$. Since $\k$ is weakly compact in $V$ there is an embedding $j:N\to M$ with critical point $\k$. I may assume that $M=\set{j(h)(\k)\st h\in N}$, by replacing $M$ with the collapse of this set if necessary. It follows, as in \FFFWC, that $M^{\ltk}\of M$. In $M$, the forcing $j(\P)$ factors as $\P*\Ptail$, where $\Ptail$ is the lottery preparation using $j(f)\restrict[\k,j(\k))$; this is $\ltk$-strategically closed in $M[G]$. Furthermore, it is not difficult to establish that $(M[G])^{\ltk}\of M[G]$. Thus, since $M[G]$ has size $\k$ we may by diagonalization construct in $V$ an $M[G]$-generic $\Gtail\of\Ptail$, and thereby lift the embedding to $j:N[G]\to M[j(G)]$ where $j(G)=G*\Gtail$. Since $\dot T_G=T$ we know $T\in N[G]$. Thus, $j(T)$ is a $j(\k)$-tree in $M[j(G)]$. Any element on the $\k^\th$ level of $j(T)$ provides a $\k$-branch through $T$ in $M[j(G)]$, and hence in $V[G]$. So $\k$ has the tree property in $V[G]$ and therefore remains weakly compact.\qed

\theorem Lottery Preparation Theorem. The lottery preparation of a measurable cardinal $\k$ satisfying $2^\k=\k\plus$ preserves the measurability of $\k$. \ref\MeasurableLottPrep

\proof This argument is similar to \FFMeasurable. Suppose $\k$ is measurable and $2^\k=\k\plus$ in $V$, and that $G$ is $V$-generic for the lottery preparation of $\k$ defined relative to some function $f$. I will show that every ultrapower embedding $j:V\to M$ by a measure on $\k$ in $V$ lifts to an embedding in $V[G]$. Given such an embedding, factor the forcing $j(\P)$ as $\P*\Ptail$, where $\Ptail$ is the lottery preparation defined in $M[G]$ using $j(f)\restrict[\k,j(\k))$. Below a condition which opts for trivial forcing in the stage $\k$ lottery, the forcing $\Ptail$ is $\ltek$-strategically closed in $M[G]$. Furthermore, standard arguments establish that $M[G]^\k\of M[G]$ in $V[G]$ and a counting argument establishes that $\card{j(\k\plus)}^V=\k\plus$. Thus, by diagonalization, we can in $V[G]$ construct an $M[G]$-generic $\Gtail\of\Ptail$; one simply lines up all the dense sets into a $\k\plus$-sequence and meets them one-by-one, following the strategy in order to get through the limit stages. Thus, in $V[G]$ the embedding lifts to $j:V[G]\to M[j(G)]$ where $j(G)=G*\Gtail$, as desired.\qed

\theorem Lottery Preparation Theorem. The lottery preparation of a strong cardinal $\k$ satisfying $2^\k=\k\plus$, defined relative to a function with the strong Menas property, preserves the strongness of $\k$.
\ref\Strong

\proof What is more, the result is completely local; following \FFFStrong, I will show that if $\k$ is $\l$-strong in $V$ then this is preserved to the lottery preparation $V[G]$. Suppose $j:V\to M$ witnesses the $\l$-strongness of $\k$ and the Menas property of $f$, so that $V_\l\of M$ and $\d=\card{V_\l}<j(f)(\k)$. By using the induced canonical extender if necessary, I may assume that $M=\set{j(h)(s)\st h\in V\and
s\in\d^{<\w}}$. Since $\d<j(f)(\k)$, below the condition $p\in j(\P)$ which opts for trivial forcing in the stage $\k$ lottery the forcing factors as $\P*\Ptail$, where $\Ptail$ is $\lted$-strategically closed. Let $\X=\set{j(h)(\k,\d)\st h\in V}\elesub M$ be the seed hull of $\<\k,\d>$, and factor the embedding as $$\trianglediagram{V}{j_0}{j}{M_0}{k}{M}$$ where
$k:M_0\to M$ is the inverse of the collapse of $\X$. Since $\k$ is in $\X$, it follows that $p\in\X$, and so $k(p_0)=p$ for some $p_0\in j_0(\P)$. Similarly, since $\d$ is in $\X$, we know that $k(\d_0)=\d$ for some $\d_0<j_0(\k)$. Also, since $\k$ is in $\X$ we know that $\cp(k)>\k$.  The embedding $j_0:V\to M_0$, being generated by the seed $\<\k,\d_0>$, is simply an ultrapower by a measure on $\k$, and therefore lifts by the diagonalization argument of \FFMeasurable\ to an embedding $j_0:V[G]\to M_0[j_0(G)]$ below the condition $p_0$. It must be that $j_0(G)$ factors as $G*\Gtail^{M_0}$, where $\Gtail^{M_0}\of\Ptail^{M_0}$ is $M_0[G]$-generic for ${\lte}\d_0$-strategically closed forcing. 

It remains to lift the embedding $k$ to $M_0[j_0(G)]$. Since $cp(k)>\k$, certainly $k$ lifts to $k:M_0[G]\to M[G]$. For the rest it suffices to show that $k\image \Gtail^{M_0}\of\Ptail$ is $M[G]$-generic. But every open dense set $D\in M[G]$ for $\Ptail$ has the form $k(\vec D)(s)$ for some $\vec D=\<D_\s\st \s\in\d_0^{<\w}>$ in $M_0[G]$ and $s\in \d^{<\w}$, where each $D_\s$ is an open dense subset of $\Ptail^{M_0}$. Since $\Ptail^{M_0}$ is ${\lte}\d_0$-distributive in $M_0[G]$, the intersection $\Dbar=\intersect_\s \vec D_\s\in M_0[G]$ remains open and dense. Furthermore, $k(\Dbar)\of D$. Thus, since $\Gtail^{M_0}$ is $M_0[G]$-generic, $k\image \Gtail^{M_0}$ meets $D$, as desired. Consequently, $k$ lifts fully to $k:M_0[j_0(G)]\to M[k(j_0(G))]$, where $k(j_0(G))$ is the filter generated by $k\image j_0(G)$. The composition $k\circ j_0$ provides a lift of $j$ to $j:V[G]\to M[j(G)]$. Since $V_\l\of M$ and $G\in M[j(G)]$, it follows that $(V[G])_\l\of M[j(G)]$, and so $\k$ is still $\l$-strong in $V[G]$.\qed

\theorem Lottery Preparation Theorem. The lottery preparation of a strongly compact cardinal $\k$, defined relative to a function with the Menas property, preserves the strong compactness of $\k$. \ref\STRCLottPrep

\proof What is more, following \FFFSTRC\ I will show that every strong compactness measure in the ground model extends to a measure in the forcing extension. Suppose that $f$ has the Menas property in $V$ (e.g. perhaps $f$ was added by fast function forcing over a smaller model), that $G\of\P$ is $V$-generic for the lottery preparation relative to $f$, that $\l\geq\k$ and that $\mu_0$ is a fine measure on $P_\k\l$ in $V$. Let $\theta\geq 2^{\l^{<\k}}$, and pick $j:V\to M$ a $\theta$-strongly compact embedding, the ultrapower by a fine measure $\eta$ on $P_\k\theta$ in $V$ witnessing the Menas property on $f$. As in \FFFSTRC, fix a seed $s_0$ for $\mu_0$, so that $X\in\mu_0\iff s_0\in j(X)$ for $X\of P_\k\l$ in $V$. In particular, $j\image\l\of s$ by fineness. Let $s=[\id]_\eta$, and $\d=\card{s}^M$. Thus, $j\image\theta\of s\in j(P_\k\theta)$ and $\theta\leq\d<j(\k)$. Now, in $j(\P)$, let $p$ be the condition which opts in the stage $\k$ lottery for a trivial poset. By the Menas property we know $j(f)(\k)>\d$, and so the next nontrivial stage of forcing lies beyond $\d$. In particular, below the condition $p$ the forcing factors as $\P*\Ptail$, where $\Ptail$ is $\lted$-strategically closed in $M[G]$. Force below $p$ to add $\Gtail\of\Ptail$ generically over $V[G]$, and in $V[G][\Gtail]$ lift the embedding to $j:V[G]\to M[j(G)]$, where $j(G)=G*\Gtail$. Let $\mu_0^*$ be the measure germinated by the seed $s_0$ via the lifted embedding, so that $X\in\mu_0^*\iff s_0\in j(X)$ for $X\in V[G]$. It is easy to see that $\mu_0^*$ measures every set in $V[G]$, that it extends $\mu_0$, and, since $j\image\l\of s_0$, that is fine. It remains only for me to show that $\mu_0^*\in V[G]$. As in \FFFSTRC\ there are again only $2^{\l^\ltk}$ many nice names in $V$ for subsets of $P_\k\l$ in $V[G]$, and we may enumerate them $u=\<\dot X_\a\st\a<\theta>$. Thus, $j(u)\in M$, and consequently also $j(u)\restrict s=\<\dot Y_\b\st \b\in s>\in M$, with $\Ydot_{j(\a)}=j(\Xdot_\a)$ for $\a<\theta$. Let $t=\set{\b\in s\st s_0\in (\Ydot_\b)_{j(G)}}$.  Thus, $t\of s$ and
$t\in M[j(G)]$.  Since $s$ has size $\d$ and $\Ptail$ is $\lted$-strategically closed in $M[G]$, it follows that $t\in M[G]$,
and therefore $t\in V[G]$.  Since $(\Xdot_\a)_G\in \mu_0^*\iff s_0\in j((\Xdot_\a)_G)=j(\Xdot_\a)_{j(G)}=(\Ydot_{j(\a)})_{j(G)}\iff j(\a)\in t$, it follows that $\mu_0^*$ is definable in $V[G]$ from $t$ and $j\restrict\theta$, and so $\mu_0^*\in V[G]$ as desired.\qed

\theorem. The lottery preparation of a supercompact cardinal $\k$, defined relative to a function with the Menas property, preserves the supercompactness of $\k$. 

\proof If in the previous argument one takes $j$ to be a $\theta$-supercompact embedding and $s_0=j\image\l$, then it is easy to see that the resulting measure $\mu_0^*$ is normal and fine on $P_\k\l$, and so $\k$ is $\l$-supercompact in $V[G]$, as desired.\qed

The previous two theorems are global in the sense that they assume $\k$ is {\it fully} strongly compact or {\it fully} supercompact in the ground model and conclude that $\k$ remains {\it fully} strongly compact or supercompact after the lottery preparation. But it is easy to extract from the proofs the following more local facts, where we assume the lottery preparation is made relative to a function with the appropriate amount of the Menas property:

\theorem Local Version. If $\k$ is $2^{\l^{<\k}}$-strongly compact in $V$, then after the lottery preparation $\k$ remains $\l$-strongly compact. If $\k$ is $2^{\l^{<\k}}$-supercompact in $V$, then, after the lottery preparation $\k$ remains $\l$-supercompact. Indeed, if $\k$ is $2^{\l^{<\k}}$-strongly compact and $\l$-supercompact in $V$, then after the lottery preparation $\k$ remains $\l$-supercompact.

A completely local result, in which the very same large cardinal assumption made in $V$ is preserved to $V[G]$, is possible if one is willing to pay a slight \GCH\ penalty:

\theorem Completely Local Version. The lottery preparation of a $\l$-supercompact cardinal $\k$ with $2^{\l^\ltk}\!=\l\plus$, defined relative to a function with the Menas property, preserves the $\l$-supercompactness of $\k$.\ref\Diagonal

\proof This argument follows the diagonalization technique used in \LocalB. Suppose that $j:V\to M$ is a $\l$-supercompact embedding in $V$ such that $j(f)(\k)>\l$ and that $V[G]$ is the lottery preparation of $\k$. By opting for trivial forcing in the stage $\k$ lottery, we may factor the forcing $j(\P)$ as $\P*\Ptail$ where $\Ptail=\P_{\l,j(\k)}$ is $\ltel$-strategically closed in $M[G]$. Standard arguments establish that $M[G]$ is closed under $\l$-sequences in $V[G]$, and a simple counting argument shows that there are at most $2^{\l^{<\k}}=\l\plus$ many subsets of $\Ptail$ in
$M[G]$, counted in $V[G]$.  Thus, I may line up the open dense subsets of $\Ptail$ in $M[G]$ into a $\l\plus$-sequence, and construct by diagonalization a descending $\l\plus$-sequence of conditions, according to the strategy, which meets every open dense set on the list. Every initial segment of the sequence is in $M[G]$, and so the
strategy ensures that the construction can proceed through any limit stage. Thus, in $V[G]$ I construct an $M[G]$-generic $\Gtail\of\Ptail$ and thereby lift the embedding to $j:V[G]\to M[j(G)]$. So $\k$ remains $\l$-supercompact in $V[G]$, as desired.\qed

The measures which exist after the lottery preparation $V[G]$ enjoy a special relationship with the measures from the ground model. Namely, I have shown that under suitable hypothesis every supercompactness or strong compactness measure in the ground model extends to a measure in the forcing extension; amazingly, the converse also holds. 

For the following theorem, assume that the first element of the domain of the function $f$ used to define the iteration is very small, say, below the least weakly compact limit of weakly compact cardinals, and that $f(\b)\geq\b$. 

\theorem. Below a condition, the lottery preparation creates no new  measurable, strong, Woodin, strongly compact, or supercompact cardinals. In addition, it does not increase the degree of strong compactness or supercompactness of any cardinal. And except possibly for certain limit ordinals of small cofinality, it does not increase the degree of strongness of any cardinal. The reason for each of these facts is that every measure in the forcing extension which concentrates on a set in the ground model extends a measure from the ground model.
\ref\Goodnews

\proof Suppose that $G$ is generic for the lottery preparation below the condition which opts, in the very first lottery at stage $\b$, to add a Cohen subset to $\b$. If $\g$ is the next element of the domain beyond $f(\b)$, then the forcing factors as $\add(\b,1)*\P_{\g,\k}$, where $\P_{\g,\k}$ is the remainder of the preparation. Thus, this is forcing of size $\b$ followed by forcing which is $\lteb\plus$-strategically closed. Such forcing is said in \cite[Ham98] to admit a {\df gap} at $\g=\b\plus$. As I explained in Remark \FFFGap, The Gap Forcing Corollary of \cite[Ham$\infty$] asserts that after such forcing every embedding $j:V[G]\to M[j(G)]$ with the property that $M[j(G)]$ is closed under $\g$-sequences---and this includes any ultrapower embedding with critical point $\k$ by a measure on any set, because such embeddings are always closed under $\k$-sequences---lifts an embedding from the ground model. That is, $M\of V$ and $j\restrict V:V\to M$ is definable in $V$. If $j$ was the ultrapower by some measure $\mu$ concentrating on a set $D\in V$, then $s=[\id]_\mu\in j(D)\in M$ is a seed for $\mu$ in the sense that $X\in\mu\iff s\in j(X)$. Consequently, $\mu\intersect V$ is definable in $V$ from $j\restrict V$. If $\mu$ is measure on $\k$, or a strong or supercompactness measure in $V[G]$, then it is not difficult to see that $\mu\intersect V$ is the corresponding kind of measure in $V$. If the original embedding $j$ was a natural $\l$-strongness extender embedding for $\l$ either a successor ordinal or a limit ordinal of cofinality above $\g$, then \cite[Ham$\infty$] shows that the restricted embedding $j\restrict V$ witnesses the $\l$-strongness of $\k$ in $V$. It follows that Woodin cardinals also cannot be created. So the theorem is proved.\qed

Thus, the lottery preparation is a gentle one; measures in the lottery preparation extension are closely related to measures in the ground model. And Remark \FFFGap\ shows that the same is true of fast function forcing, provided that it is prefaced by some small forcing. 

\section Indestructibility After the Lottery Preparation

Now I come to the main contribution of this paper, namely, that the lottery preparation makes large cardinals indestructible.  Laver's \cite[Lav78] original preparation---my inspiration, of course---showed spectacularly that any supercompact cardinal can be made highly indestructible. Gitik and Shelah \cite[GitShl89] extended the analysis to strong cardinals. Here, the lottery preparation unifies and generalizes these results by providing a uniform preparation which works with any large cardinal, whether it is supercompact, strongly compact, strong, partially supercompact, partially strongly compact, or merely measurable, and so on. In each of these cases, the lottery preparation makes the cardinal indestructible by a variety of forcing notions, depending on the strength of the cardinal in the ground model. 

Of course, the lottery preparation will make a supercompact cardinal $\k$ fully indestructible by $\ltk$-directed closed forcing, and a strong cardinal $\k$ (such that $2^\k=\k\plus$) indestructible by $\ltek$-distributive forcing, and more. These arguments, given below, essentially follow the corresponding results in \cite[Lav78] and \cite[GitShl89]. The level-by-level result, lacking in the Laver preparation because Laver functions are not available level-by-level, is nevertheless possible with the lottery preparation, which requires no Laver function. Thus, with a bit of the \GCH, even partially supercompact cardinals can be made indestructible.

The most significant contribution of this paper, however, concerns the strongly compact cardinals. Because the lifting arguments involved in Laver's theorem simply fail when used with a strongly compact embedding, it has been an open question for some time to determine the degree to which a strongly compact cardinal can be preserved by forcing. Indeed, Menas \cite[Men74] seems very concerned with this question. Apter concludes his paper \cite[Apt97] with questions asking how much indestructibility is possible with a strongly compact non-supercompact cardinal. Until now, there have been no nontrivial instances of an arbitrary strongly compact cardinal being preserved by forcing. Progress had been made in the special case of a strongly compact limit of supercompact cardinals (for which a Menas function always exists): Apter \cite[Apt96], using Menas's technique, showed how to make such cardinals $\k$ indestructible by $\add(\k,1)$; Menas \cite[Men74] himself also seems also very close to proving this. Apter \cite[Apt97] shows how to make any strongly compact limit of supercompact cardinals indestructible by any $\ltk$-directed closed forcing which does not add a subset to $\k$.  Recently, Apter and Gitik \cite[AptGit97] proved, impressively, that any supercompact cardinal can be made into strongly compact cardinal which is simultaneously the least measurable cardinal and fully indestructible by $\ltk$-directed closed forcing. Here, I aspire to eliminate such supercompactness assumptions and work with an arbitrary strongly compact cardinal. How much indestructibility is possible? 

For the remainder of this section, assume that $V[G]$ is the lottery preparation of $\k$ relative to a function $f$ with the suitable Menas property in $V$. For example, perhaps $f$ was previously added by fast function forcing. Occasionally, for a few of the theorems, I will make the additional assumption that in fact $f$ was added in this way. 

\theorem Indestructibility Theorem. After the lottery preparation, a strongly compact cardinal $\k$ becomes indestructible by Cohen forcing $\add(\k,1)$. \ref\Cohen

\proof Suppose that $V[G]$ is the lottery preparation of $\k$ relative to the function $f$ with the strongly compact Menas property, and that $g\of\k$ is $V[G]$-generic for $\add(\k,1)^{V[G]}$. I will show that every fine measure $\mu_0$ on $P_\k\l$ in $V$ extends to a measure in $V[G][g]$.  Let $\theta\geq2^{\l^{<\k}}$, and suppose $j:V\to M$ witnesses the Menas property of $f$, so that it is the ultrapower by a fine measure $\eta$ on $P_\k\theta$ in $V$ and  $j(f)(\k)>\d=\card{s}$ where $s=[\id]_\eta$. Below the condition $p$ which opts in the stage $\k$ lottery for $\add(\k,1)^{M[G]}=\add(\k,1)^{V[G]}$, the forcing $j(\P)$ factors as $\P*\add(\k,1)*\Ptail$, where the forcing $\Ptail$ is $\lted$-strategically closed in $M[G][g]$. Force to add $\Gtail\of\Ptail$ over $V[G][g]$, and lift the embedding to $j:V[G]\to M[j(G)]$, where $j(G)=G*g*\Gtail$. Now consider the forcing $j(\add(\k,1)^{V[G]})=\add(j(\k),1)^{M[j(G)]}$.  The set $g\of\k$ is a condition in this poset, so I can force below it to add a generic $\tilde g\of j(\k)$. Since $\tilde g$ pulls back to $g$ via $j$, that is to say, since $g$ is a master condition, the embedding lifts in $V[G][g][\Gtail][\tilde g]$ to $j:V[G][g]\to M[j(G)][j(g)]$ where $j(g)=\tilde g$. Select a seed $s_0$ for $\mu_0$ as in Theorem \FFFSTRC, and let $\mu_0^*$ be the $V[G][g]$-measure germinated by $s_0$ via the lifted embedding, so that $X\in\mu_0^*\iff s_0\in j(X)$. I want to show that $\mu_0^*$ is in $V[G][g]$. Enumerate $u=\<\Xdot_\a\st\a<\theta>$ the names for subsets of $P_\k\l$ in $V[G][g]$. Again, let $t$ be the set of all $\b\in s$ such that $s_0\in (\Ydot_\b)_{j(G)*j(g)}$, where $j(u)\restrict s=\<\Ydot_\b\st\b\in s>$.  Since $t\of s$ and the tail forcing is $\lted$-strategically closed, it follows that $t\in M[G][g]$, and consequently $t\in V[G][g]$. Using $t$ and $j\restrict\theta$ as in Theorem \STRCLottPrep, I conclude that $\mu_0^*$ is in $V[G][g]$, as desired.\qed

Let me introduce now another kind of forcing for which strongly compact cardinals will become indestructible.  For any set $S$ which is in a normal measure on $\k$, the {\df club forcing} $\Q_S$ will add a club $C\of\k$ such that $C\intersect\inacc\of S$; conditions are closed bounded sets $c\of\k$ such that $c\intersect\inacc\of S$, ordered by
end-extension.  For every $\b<\k$, the set of such $c$ which mention an element above $\b$ is a $\lteb$-closed open dense set, since one can simply take the union of a $\b$-chain of such conditions and add the supremum to obtain a stronger condition; the supremum cannot be inaccessible since it is above $\b$ but was reached by the
$\b$-sequence. Thus $\Q_S$ preserves all cardinals and cofinalities. A variant, the {\df coherent club forcing} $\Qtilde_S$, is meant directly to follow a lottery preparation or other iteration, and imposes the additional requirement that whenever $\d$ is an inaccessible cluster point of $C$, then the preceding iteration added the club
$C_\d=C\intersect\d$ by forcing with $\Q_{S\intersect \d}$ at stage $\d$.  This forcing adds a coherent system of clubs which reflect at their inaccessible cluster points.

In the next theorem I will need the simple fact (proved also in \cite[Men74]) that if $\mu$ is a normal fine measure on $P_\k\l$ and $\eta$ is a fine measure on $P_\k\theta$ for some $\theta\geq\l$ of cofinality at least $\k$, then the product measure $\mu\cross\eta$ is isomorphic to a fine measure on $P_\k\theta$, and the resulting $\theta$-strongly compact embedding $j:V\to M$ is closed under $\l$-sequences; in particular, $j\image\l\in M$. To see why this is true, consider the commutative diagram corresponding to the product measure $\mu\cross\eta$:
$$\trianglediagram{V}{j_\mu}{j}{M_\mu}{k}{M}$$ where $k$ is the ultrapower of $M_\mu$ by $j_\mu(\eta)$. Every element of $M_\mu$ has the form $j_\mu(h)(j_\mu\image\l)$ for some $h\in V$, and every element of $M$ has the form $k(F)(s)$, where $F\in M_\mu$ and $s=[\id]_{j_\mu(\eta)}$. Thus, since $k(j_\mu\image\l)=j\image\l$,
every element of $M$ has the form $j(h)((j\image\l),s)$. Let $t$ be the element of $j(P_\k\theta)$ which is obtained in $M$ by simply placing a copy of $j\image\l$ at the top of $s$, separated by a brief gap. From $t$ one can recover both $s$ and $j\image \l$, so every element of $M$ has the form $j(h)(t)$ for some function $h\in V$. In the seed terminology of \cite[Ham97], the seed $t$ generates all of $M$.  It follows that $j$ is the ultrapower by the corresponding measure $\bar\eta$, defined by $X\in\bar\eta\iff t\in j(X)$, and consequently that $\bar\eta$ is isomorphic to $\mu\cross\eta$.  Since $j\image\theta\of s\of t$, it follows that $\bar\eta$ is a fine measure on $P_\k\theta$.  Thus, $j$ is a $\theta$-strongly compact embedding, as desired. Since $j\image\l\in M$, it follows that $M$ is closed under $\l$-sequences, since any $\l$-sequence $\<j(f_\a)(t)\st\a<\l>$ is equal to $j(F)(t,j\image\l)$ where $F(\s,\t)=\<f_\a(\s)\st\a\in\t>$, and is therefore in $M$. So I have proved the fact that I need. The argument also works to show that if $\mu$ is merely a fine measure on $P_\k\l$, then still $\mu\cross\eta$ is isomorphic to a fine measure on $P_\k\theta$.

For the purposes of the next theorem, let me say that a subset $S$ of the strongly compact cardinal $\k$ is {\df special} when for arbitrarily large $\theta$ it is in the induced normal measure of a $\theta$-strong compactness embedding witnessing the Menas property of $f$. If $2^\k=\k\plus$ and the function $f$ was added by fast function forcing, then the special sets include any set in a normal measure on $\k$ in the original ground model. To see why this is so, suppose that $S$ is in a normal measure on $\k$ in $\Vbar$ and $V=\Vbar[f]$ is the fast function extension. By \FFMeasurable\ we know that $S$ is in a normal measure in $V$, and in $V$ we can take a product of this measure with any $\theta$-strong compactness measure to get a $\theta$-strong compactness measure with $\k\in j(S)$. Then, by the Flexibility Theorem \Flexibility, we can modify $j(f)$ and factor the embedding in the manner of Theorem \STRCMenas\ to ensure that $j(f)(\k)$ is large enough to witness the Menas property; since $S\in\Vbar$, these modifications do not affect whether $\k\in j(S)$, and so the desired hypothesis is obtained. Therefore, we have numerous interesting sets $S$ which are special.

\theorem Indestructibility Theorem. After the lottery preparation, a strongly compact cardinal $\k$ becomes indestructible by the club forcing $\Q_S$, and by the coherent club forcing $\Qtilde_S$, for any special set $S$ in $V$.\ref\Club

\proof Suppose that $V[G]$ is the lottery preparation of the strongly compact cardinal $\k$ relative to $f$, that $S\in V$ is special and that $C$ is the generic club added by $\Q_S$. I want to show that $\k$ remains strongly compact in $V[G][C]$; indeed, I will show that all the strongly compact measures in $V$ extend to measures in $V[G][C]$. Fix any $\l\geq\k$ and choose $\theta\geq 2^{\l^\ltk}$ and a $\theta$-strong compactness measure $\eta$ on $P_\k\theta$ which witnesses the speciality of $S$. Thus, the corresponding embedding $j:V\to M$ has $\k\in j(S)$ and $j(f)(\k)>\d=\card{s}$ where $s=[\id]_\eta$. 

Let $p$ be the condition which opts for $\Q_S$ in the stage $\k$ lottery, so that below $p$ the forcing $j(\P)$ factors as $\P*\Q_S*\Ptail$. Add $\Gtail\of\Ptail$ generically over $V[G][C]$ (this is $\lted$-strategically closed in $M[G][C]$), and lift the embedding to $j:V[G]\to M[j(G)]$, where $j(G)=G*C*\Gtail$. Because $C$ is explicit in $j(G)$, we know $C\in M[j(G)]$; and since $\k\in j(S)$ the set $\Cbar=C\union\set{\k}$ is available as a condition in $j(\Q_S)$. In the coherent club context, in which $C$ is generic for $\Qtilde_S$, observe similarly that $\Cbar$ is a condition in $j(\Qtilde_S)$, since $\Cbar\intersect \k=C$, and this is the generic which was used at stage $\k$ in $j(G)$.  In either context, let $j(C)$ be $V[G][C][\Gtail]$-generic below $\Cbar$, and lift the embedding to $j:V[G][C]\to M[j(G)][j(C)]$ in $V[G][C][\Gtail][j(C)]$.  For any fine measure $\mu_0$ on $P_\k\l$ in $V$, we can find by the cover property for $j$ an element $s_0$ such that $X\in\mu_0\iff s_0\in j(X)$ for $X\of P_\k\l$ in $V$.  Use this same seed to germinate a measure with respect to the lifted embedding $j$ according to the rule $X\in \mu_0^*\iff s_0\in j(X)$ for $X\of P_\k\l$ in $V[G][C]$. It is clear that $\mu_0^*$ extends $\mu_0$, and the argument of Theorem \STRCLottPrep, using the enumeration $u$ of the names for subsets of $P_\k\l$, shows that $\mu_0^*$ is in $V[G][C]$, as desired.\qed

I would like now to amuse the reader by applying the idea of the previous theorem with the set $S$ of non-measurable cardinals. 

\corollary. Assume that the \GCH\ holds. Then, while preserving the strong compactness of any strongly compact cardinal $\k$, one can add a club $C\of\k$ which contains no measurable cardinals. Furthermore, this can be done while preserving all cardinals and cofinalities, and while neither creating nor destroying any measurable cardinals.

\proof We first add a fast function $f\from\k\to\k$ over $V$. This preserves all cardinals and cofinalities and neither creates nor destroys any measurable cardinals. Next, we will force over the modified lottery preparation $\P$, in which forcing is allowed in the stage $\g$ lottery only when, in addition to the earlier requirement that it is in $H(f(\g)\plus)$ and for every $\b<\g$ it is $\ltb$-strategically closed, but also that it preserves all cardinals and cofinalities and does not destroy any measurable cardinals. Suppose that $G\of\P$ is $V[f]$-generic for this modified preparation. Since we can arrange the preparation to admit a very low gap, by Remark \FFFGap\ the forcing does not 
create any measurable cardinals. By the remarks preceding the previous theorem, the set $S$ of non-measurable cardinals in $V$ is special in $V[f]$, and consequently, since the modifications to the lottery preparation here do not
create any difficulties in the lifting argument of the previous theorem, it follows that $\k$ remains strongly compact in 
in $V[f][G][C]$, where $C\of\k$ is the club added by forcing with $\Q_S$ over $V[f][G]$. Since $C\intersect\inacc\of S$, it follows that $C$ contains no cardinals which are measurable in $V$. Since no measurable cardinals are created, $C$ contains no cardinals which are measurable in $V[f][G]$ or in $V[f][G][C]$. Because the forcing at every stage preserves cardinals and cofinalities, the standard reverse Easton iteration arguments establish that the entire iteration also preserves cardinals and cofinalities. Finally, I will show that all measurable cardinals are preserved. Certainly the measurable cardinals in $V$ above $\k$ are preserved. 
Also, I have argued that the strong compactness of $\k$ itself is preserved. So suppose that $\g<\k$ is a measurable 
cardinal in $V$, and hence also $V[f]$. Since the forcing after stage $\g$ in $\P$ is strategically closed up to the next inaccessible cardinal, it cannot affect the measurability of $\g$. Also, forcing at stage $\g$ is only allowed when it preserves the measurability of $\g$. Thus, it suffices to show that $\g$ is measurable in $V[f][G_\g]$. There are two cases. First, it may happen that $f\image\g\of\g$. In this case, the forcing up to stage $\g$ is exactly the modified lottery preparation of $\g$, which by the argument of \MeasurableLottPrep\ preserves the measurability of $\g$. Second, alternatively, it may happen that for some $\b<\g$ we have $f(\b)\geq\g$. So there is no forcing between stage $\b$ and $\g$.  In this case, the forcing up to stage $\b$ is small relative to $\g$, and therefore preserves the measurability of $\g$, and the forcing at stage $\b$ was only allowed provided that it also preserved the measurability of $\g$, so $\g$ is measurable in $V[f][G_\g]$, as desired.\qed

In the previous argument, if one uses coherent club forcing one obtains also a whole sequence of clubs $C_\g\of\g$ for lots of $\g\leq\k$, all disjoint from the measurable cardinals, with the coherency property, so that whenever $\b$ is an inaccessible cluster point of $C_\g$, then $C_\b$ exists and $C_\g\intersect\b=C_\b$. 

Let me introduce another forcing notion for which strongly compact cardinals become indestructible. The {\df long Priky forcing} poset $\Q_F$, where $F$ is a $\k$-complete filter on $\k$, consists of conditions $\<s,A>$, where $s\in [\k]^\ltk$ and $A\in F$, ordered in the Prikry manner, so that $\<s,A>\leq\<t,B>$ when $s$ end-extends $t$, $A\of B$, and $s\setminus t\of B$. This forcing adds a single set $g\of\k$ such that every set in $F$ contains a tail of $g$. It is $\ltk$-directed closed and has the $\k\plus$-chain condition; so all cardinals and cofinalities are preserved. Define that a set $z$ is {\df accessible} to an embedding $j:\Vbar\to\Mbar$ when $z\in\Mbar$. 

\theorem Indestructibility Theorem. After the lottery preparation relative to a fast function, a strongly compact cardinal $\k$ becomes indestructible by long Prikry forcing $\Q_F$ for any $\k$-complete filter $F$ on $\k$ which is accessible to a strongly compact embedding.

\proof Suppose $V[f][G]$ is the lottery preparation of $\k$ obtained by first adding the fast function $f$ and $g\of\Q_F$ is $V[f][G]$-generic for long Prikry forcing $\Q_F$ on the $\k$-complete filter $F$ which is accessible to a strongly compact embedding in $V[f][G]$. By the product measure argument preceding Theorem \Club, it follows by taking products that $F$ is accessible to $\theta$-strongly compact embeddings for all sufficiently large $\theta$. Suppose $\theta\geq 2^{\l^\ltk}$ and $F$ is accessible to the ultrapower embedding $j:V[f][G]\to M[j(f)][j(G)]$ by the $\theta$-strong compactness measure $\eta$. By the closure considerations on the forcing, $F\in M[f][G][\tilde g]$, where $\tilde g\of\R$ is the (possibly trivial) stage $\k$ forcing in $j(G)$. Let $\d=\card{s}$ where $s=[\id]_\eta$. By Lemma \STRCMild, we may replace $\eta$ if necessary with an isomorphic measure and assume $s\in M$. By the Flexibility Theorem applied to $j\restrict V[f]$, there is another embedding $j^*:V[f]\to M[j^*(f)]\of M[j(f)]$ such that $j^*(f)(\k)>\d$. Since $\R$ was allowed at stage $\k$, it is in $H(j(f)(\k)\plus)^{M[j(f)][G]}$, and consequently, by the closure of $\ftail$, it is in $M[f][G]$. And since also $\Q_F\in M[f][G][\tilde g]$ and the poset $\R*\Q_F$ satisfies the required closure conditions, it is allowed to appear in the stage $\k$ lottery of $j^*(\P)$. Let $p$ be the condition which opts for $\R*\Q_F$ at stage $\k$, so that below $p$ the forcing
$j^*(\P)$ factors as $\P*(\R*\Q_F)*\Ptail$, where $\Ptail$ is $\lted$-strategically closed in $M[j^*(f)][G][\tilde g][g]$. Force 
to add $\Gtail^*\of\Ptail$ over $V[f][G][g]$ and lift the embedding to $j^*:V[f][G]\to M[j^*(f)][j^*(G)]$ where $j^*(G)=G*(\tilde g*g)*\Gtail^*$. Now consider the $j^*(\Q_F)$ forcing. Enumerate $F=\<X_\a\st\a<2^\k>$. Since $j\image\theta\of s$, it follows that $j^*(F)\restrict s=\<Y_\b\st \b\in s\restrict j(2^\k)>$  provides a cover of $j^*\image F$ of size at most $\d$. Since the sets $Y_\b$ are all in $j^*(F)$, and $j^*(F)$ is a $j^*(\k)$-complete filter in $M[j^*(f)][j^*(G)]$, I can intersect them all to obtain a set $Y=\intersect\set{Y_\b\st\b\in s\restrict j(2^\k)}\in j^*(F)$. Since $j^*(X_\a)=Y_{j^*(\a)}$, it follows that $Y\of j^*(X)$ for any $X\in F$. Thus, $\<g,Y>$ is a condition in $j^*(\Q_F)$ with the property that any $\<t,A>\in g$ has $\<g,Y>\leq\<t,j^*(A)>=j^*(\<t,A>)$; that is, it is a master condition. Force below it to add the generic $j^*(g)$, and in $V[f][G][g][\Gtail^*][j^*(g)]$ lift the embedding to $j^*:V[f][G][g]\to M[j(f)][j^*(G)][j^*(g)]$. If $\mu_0$ is any fine measure on $P_\k\l$ in $V$, then we may find a seed $s_0\in M$ for $\mu_0$ as in \STRCLottPrep, and let $\mu_0^*$ be the measure germinated by $s_0$ via $j^*$. Certainly $\mu_0^*$ extends $\mu_0$ and the argument of \STRCLottPrep\ involving the enumeration $u$ of the nice names in $V$ for subsets of $P_\k\l$ in $V[f][G][g]$ shows that $\mu_0^*$ lies in $V[f][G][g]$, as desired.\qed

In the case that a strongly compact cardinal has some nontrivial degree of supercompactness, this partial supercompactness can be used to obtain more indestructibility for the full strong compactness. 

\theorem Indestructibility Theorem. If $\k$ is strongly compact and $\l$-supercompact in the ground model, then after the lottery preparation relative to a fast function, both of these properties are indestructible by any $\ltk$-directed closed forcing of size less than or equal to $\l$.\ref\STRCSC

\proof First let me show the preliminary claim in $V[f]$ that for any $\theta$ there is a $\l$-closed $\theta$-strongly compact embedding witnessing the Menas property of $f$. Begin with the argument preceding Theorem \Club, which produces a $\theta$-strongly compact embedding $j:V[f]\to M[j(f)]$, the ultrapower by a fine measure $\eta$ on $P_\k\theta$, which is closed under $\l$-sequences. In particular, $j\image\l\in M[j(f)]$. By Remark \FFFGap, it must be that $j\image\l\in M$. By Lemma \STRCMild, we may assume that $s=[\id]_\eta\in M$, and moreover that $j\image\l$ is directly coded into the top elements of $s$. Now, by the Flexibility Theorem \Flexibility, there is another embedding $j^*:V[f]\to M[j^*(f)]\of M[j(f)]$ with $j^*(f)(\k)>\card{s}$. Let $\X=\set{j^*(h)(s)\st h\in V[f]}\elesub M[j^*(f)]$ be the seed hull of $s$ with respect to $j^*$, and $j_0:V[f]\to M_0[j_0(f)]$ the induced factor embedding, with $j_0=\pi\compose j$ where $\pi$ is the Mostowski collapse of $\X$. It follows that $s_0=\pi(s)$ generates all of $M_0[j_0(f)]$, and so $j_0$ is the ultrapower by the $\theta$-strongly compact measure $\eta_0$ germinated by $s$ via $j$ (or by $s_0$ via $j_0$). Since $j\image\l\in\X$, it follows that $j_0\image\l\in M_0[j_0(f)]$. Furthermore, $j_0(f)(\k)>\card{s_0}$. In particular, $j_0$ is a $\l$-closed $\theta$-strongly compact embedding which witnesses the Menas property of $f$, so the preliminary claim is proved. 

Continuing with the main argument now, fix any $\bar\l\geq\l$ and any $\theta\geq 2^{2^{\bar\l^{{<}\k}}}$, and suppose that $j:V[f]\to M[j(f)]$ is a $\l$-closed $\theta$-strongly compact embedding witnessing the Menas property of $f$. We may assume $j\image\l\in M$ and $j(f)(\k)>\d=\card{s}$ where $s=[\id]_\eta\in M$. Suppose $g\of\Q$ is $V[f][G]$-generic for the $\ltk$-directed closed forcing $\Q$ of size at most $\l$. By the techniques used previously I can lift the embedding to $j:V[f][G]\to M[j(f)][j(G)]$ in $V[f][G][g][\Gtail]$ such that the generic $j(G)=G*g*\Gtail$ opts for the forcing $\Q$ at stage $\k$, and the next stage of forcing is beyond $\d$. I know that $\Q\in M[j(f)][G]$ since $M[j(f)]$ is closed under $\l$-sequences and $\Q$ has size at most 
$\l$. It is allowed to appear in the stage $\k$ lottery because $j(f)(\k)>\d\geq\theta$. Consider now the forcing $j(\Q)$. 
From $g$ and $j\image\l$ we can construct $j\image g$ in $M[j(f)][j(G)]$. And since this set has size $\l<j(\k)$ and is directed, there is a condition $p\in j(\Q)$ below every element of $j\image g$. Force below $p$ to add $j(g)$, and lift the embedding to $j:V[f][G][g]\to M[j(f)][j(G)][j(g)]$ in $V[f][G][g][\Gtail][j(g)]$. Let $\mu$ be the set of all $X\of P_\k\l$ in $V[f][G][g]$ such that $j\image\l\in j(X)$. It is easy to see that this is normal and fine. Furthermore, the arguments of the previous theorems show that $\mu\in V[f][G][g]$. Consequently, $\k$ is $\l$-supercompact there. Finally, any fine measure $\mu_0$ on $P_\k\bar\l$ extends to a fine measure $\mu_0^*$ in $V[f][G][g]$ by the arguments given previously, so $\k$ is strongly compact there as well, as desired. Indeed, I have shown that every $\l$-supercompactness measure in $V$ extends to a measure in $V[f][G][g]$, and every strong compactness measure in $V$ extends to a measure in $V[f][G][g]$.\qed

Notice that while the lottery preparation uses strategically closed forcing at every stage, I only claim preservation by $\ltk$-directed closed forcing in the previous theorem. This cannot be generalized to include all $\ltk$-strategically closed forcing, because the forcing which adds a non-reflecting stationary subset to $\k$ is $\ltk$-strategically
closed, but always destroys even the weak compactness of $\k$. The reason for using strategically closed forcing in the lottery preparation is to allow for the forcing such as $\Q_S$, which is not generally $\ltk$-closed, while simultaneously retaining the distributivity of the tail forcing $\Ptail$.

In any case, it follows from the previous argument that for a supercompact cardinal, the lottery preparation accomplishes everything that the original Laver preparation was meant to accomplish:

\corollary. After the lottery preparation, a supercompact cardinal $\k$ becomes indestructible by any $\ltk$-directed closed forcing. 

\proof If one uses a fast function, this corollary follows immediately from the previous theorem. Let me illustrate, nevertheless, how one can directly follow Laver's \cite[Lav78] original argument in the lottery context, while assuming only the Menas property on $f$. Suppose $g\of\Q$ is $V[G]$-generic for some $\ltk$-directed closed forcing $\Q$. Fix any $\l$ and let $j:V\to M$ be a $\theta$-supercompact embedding for some $\theta\geq 2^{\l^\ltk}, \card{\Q}$ which witnesses the Menas property of $f$, so that $j(f)(\k)>\theta$. Below a condition which opts for $\Q$ in the stage $\k$ lottery, the forcing $j(\P)$ factors as $\P*\Q*\Ptail$, where $\Ptail$ is $\ltet$-strategically closed in $M[G][g]$. Force to add a generic $\Gtail\of\Ptail$, and in $V[G][g][\Gtail]$ lift the embedding to $j:V[G]\to M[j(G)]$ where $j(G)=G*g*\Gtail$. Now, using $j\image\theta$, it follows that $j\image g\in M[j(G)]$, and so by the directed closure of $j(\Q)$ there is a master condition $q\in j(\Q)$ below every element of $j\image g$. Force to add $j(g)$ below $q$ and lift the embedding in $V[G][g][\Gtail][j(g)]$ to $j:V[G][g]\to M[j(G)][j(g)]$. Let $\mu$ be the normal fine measure on $P_\k\l$ germinated via $j$ by $j\image\l$. Since $\mu$ measures every set in $V[G][g]$, it suffices to show that $\mu\in V[G][g]$. Certainly $\mu$ is in $V[G][g][\Gtail][j(g)]$. Since the forcing $\Gtail*j(g)$ was $\ltet$-directed closed in $V[G][g]$, it could not have added $\mu$. So $\mu$ is in $V[G][g]$, as desired.\qed

\theorem Improved Version. After the lottery preparation relative to a fast function and any further $\ltk$-directed closed forcing, every supercompactness measure on $\k$ from the ground model extends to a measure in the forcing extension, and every supercompactness measure in the forcing extension extends a measure from the ground model. Furthermore, if the \GCH\ holds, then every sufficiently large supercompactness embedding from the ground model lifts to the extension.

This is an improvement over the Laver preparation, through which one can lift an embedding $j:V\to M$ only when $j(\ell)(\k)$ is appropriate. 

\proof The first half of the first sentence follows immediately from the proof of \STRCSC. The second half follows from Remark \FFFGap. The second sentence follows by the diagonalization technique of Theorems \Diagonal\ and \LocalB. Specifically, if $V[f][G]$ is the lottery preparation relative to the fast function, $g\of\Q$ is $\ltk$-directed closed forcing of size at most $\l$ and $2^{\l^{<\k}}=\l\plus$, then any $\l$-supercompactness embedding $j:V\to M$ lifts to the extension $j:V[f][G][g]\to M[j(f)][j(G)][j(g)]$.\qed

The previous argument admits a completely local analogue in a way that Laver's original preparation does not. In general, one cannot perform the Laver preparation of a $\l$-supercompact cardinal $\k$ unless one has $\l$-supercompactness Laver function; but Laver's proof that such a function exists requires that $\k$ is $2^{\l^\ltk}\!$-supercompact in the ground model. Thus, it has been open whether any partially supercompact cardinal can be made indestructible, even assuming the \GCH. This question is answered by the following theorem.

\theorem Level-by-level Preparation. If $\k$ is $\l$-supercompact in $V$ and $2^{\l^\ltk}\!=\l\plus$, then after the lottery preparation the $\l$-supercompactness of $\k$ is fully indestructible by $\ltk$-directed closed forcing of size at most $\l$.\ref\SCLevel

\proof This is essentially what I actually argued in the previous theorem. To support the diagonalization argument, one only needs the Menas property on $f$.\qed

I have shown by the previous theorems that the lottery preparation makes any strongly compact cardinal $\k$ partially indestructible; but perhaps there is much more indestructibility than I have identified, so it is natural to ask:

\question. For which other natural forcing notions does a strongly compact cardinal $\k$ become indestructible after the lottery preparation? \ref\Which

Let me consider now the lottery preparation of a strong cardinal $\k$. Recall that $\k$ is {\df strong} when for every $\l$ the cardinal $\k$ is $\l$-strong, so that there is an embedding $j:V\to M$ with critical point $\k$ such that $V_\l\of M$. Gitik and Shelah \cite[GitShl89], using Woodin's \cite[CumWdn] technique for preserving a strong cardinal, showed how to make any strong cardinal indestructible by ${\lte}\k$-directed closed forcing (indeed, they improve this to weakly
$\k$-closed posets with the Prikry property). I would like to show that such indestructibility is also achieved by the Lottery preparation.

\theorem. After the lottery preparation of a strong cardinal $\k$ such that $2^\k=\k\plus$, the strongness of $\k$ becomes indestructible by ${\lte}\k$-strategically closed forcing.

\proof This is similar to the proof of Theorem \Strong, except that I will opt for the appropriate forcing at stage $\k$ in $j(\P)$. Suppose that $V[G]$ is the lottery preparation defined relative to a function $f$ with the Menas property for strong cardinals in $V$. The result is local in that only the $\l$-strongness of $\k$ in $V$ is needed to know that the $\l$-strongness of $\k$ is indestructible over $V[G]$ by any forcing notion of rank below $\l$.  Suppose $H\of\Q\in (V[G])_\l$ is generic over $V[G]$ for ${\lte}\k$-strategically closed forcing $\Q$. Fix a $\l$-strong embedding $j:V\to M$ from the ground model such that $j(f)(\k)>\l$. As in Theorem \Strong, I may assume that $M=\set{j(h)(s)\st h\in V\and s\in\d^{<\w}}$, where $\d=\beth_\l<j(\k)$. Let $p$ be the condition in $j(\P)$ which opts to force with $\Q$ in the stage $\k$ lottery. Since the next inaccessible beyond beyond $\l$ must be beyond $\d$, the forcing $j(\P)$ factors below $p$ as $\P*\Q*\Ptail$, where $\Ptail$ is $\lted$-strategically closed in $M[G][H]$. There must be some ordinal $\b<\d$ and function $h$ such that $\dot\Q=j(h)(\b)$ for some name $\dot\Q$ for $\Q$.  Let $$\X=\set{(\dot z)_{G*H}\st \dot z=j(g)(\k,\b,\d)\hbox{ for some function }g\in V}.$$ It is not difficult to verify the Tarski-Vaught criterion, so that $\X\elesub M[G][H]$.  Also, $\X$ is closed under $\k$-sequences in 
$V[G][H]$.  Note that $\k$, $\b$, $\d$, $p$, $\Q$ and $\Ptail$ are all in $\X$. Furthermore, since $\Ptail$ is $j(\k)$-c.c. in $M[G][H]$ and there are only $2^\k=\k\plus$ many functions $g:\k\to V_\k$ in $V$, there are at most $\k\plus$ many open dense subsets of $\Ptail$ in $\X$. Since $\Ptail$ is ${\lte}\k$-closed, one can perform the diagonalization argument to construct in $V[G][H]$ a filter $\Gtail\of\Ptail$ which is $\X$-generic.  Let me argue now that $\Gtail$ is also $M[G][H]$-generic. If $D$ is an open dense subset of $\Ptail$ in $M[G][H]$, then $D=\dot D_{G*H}$ for some name $\dot D\in M$. Consequently, $\dot D=j(g)(\k,\k_1,\ldots,\k_n)$ for some function $g\in V$ and $\k<\k_1<\cdots<\k_n<\d$. Let $\Dbar $ be the intersection
of all $j(g)(\k,\a_1,\ldots,\a_n)$ where $\k<\a_1<\cdots<\a_n<\d$ such that this is an open dense set in $\Ptail$. Since $\Ptail$ is ${\leq}\d$-distributive by the choice of $p$, it follows that $\Dbar $ remains open and dense. Furthermore, $\Dbar \in\X$. Thus, since $\Gtail$ is $\X$-generic, it meets $\Dbar $, and since $\Dbar \of D$, it must also meet $D$. So $\Gtail$ is $M[G][H]$-generic, and the embedding lifts to $j:V[G]\to M[j(G)]$ where $j(G)=G*H*\Gtail$. To lift the embedding through the $\Q$ forcing, it suffices to argue that $j\image H\of j(\Q)$ is $M[j(G)]$ generic. Given an open dense set $D\of j(\Q)$ in that model, $D=\dot D_{j(G)}$ for some name $\dot D=j(\vec D)(s)$ for some function $\vec D=\<D_\s\st\s\in \k^{<\w}>$ and some $s\in\d^{<\w}$. I may assume that every $D_\s$ is open and dense in $\Q$. Since $\Q$ is ${\leq}\k$-distributive, it follows that $D'=\intersect_\s D_\s$ is also open dense in $\Q$. Since $j(D')\of D$, and $H$ meets $D'$, it follows that $j\image H$ meets $D$, as desired. Thus, the embedding lifts fully to $j:V[G][H]\to M[j(G)][j(H)]$, and this lifted embedding remains $\l$-strong because $V_\l\of M$ and so
$(V[G][H])_\l\of M[G][H]\of M[j(G)][j(H)]$.\qed

One can also mimic the master condition arguments from above to obtain:

\theorem. After the lottery preparation of a strong cardinal $\k$ relative to a fast function, the strongness of $\k$ is indestructible by $\add(\k,1)$ and by $\Q_S$ and $\Qtilde_S$ whenever $\k\in j(S)$ for arbitrarily large $\l$-strong embeddings $j$. 

\section Impossibility Theorem

One might hope to generalize the previous theorems by proving that the lottery preparation or some other alternative to the Laver preparation can make any strongly compact cardinal fully indestructible. But this hope will not be fulfilled; the sad fact which I will now prove is that no preparation which naively resembles the Laver preparation can make a strongly compact non-supercompact cardinal fully indestructible.

\theorem Impossibility Theorem. The lottery preparation will always fail to make a strongly compact cardinal fully indestructible unless it was originally supercompact. In fact, any forcing which resembles the Laver preparation---an iteration of strategically closed forcing in which the next nontrivial stage of forcing lies beyond the size of the previous
one---will fail to make a strongly compact non-supercompact cardinal fully indestructible. Indeed, after adding a single Cohen real, there is no ${\lte}\omega_1$-strategically closed preparatory forcing which makes a strongly compact non-supercompact cardinal $\k$ fully indestructible.
\ref\ImpossibilityTheorem

This theorem relies on my recent work in \cite[Ham$\infty$] and \cite[Ham98], in which, as I mentioned in Remark \FFFGap, I defined that a notion of forcing admits a gap below $\k$ when it factors as $\P_1*\P_2$ where, for some $\d<\k$, $\card{\P_1}<\d$ and $\forces\P_2$ is $\lted$-strategically closed. Any kind of Laver preparation, obtained by iterating the closed forcing provided by some kind of Laver function, admits numerous gaps below $\k$. The lottery preparation admits a gap between any two lottery stages. The Impossibility Theorem \ImpossibilityTheorem, therefore, is an immediate consequence of the following theorem.

\theorem. Forcing which admits a gap below a strongly compact cardinal $\k$ cannot make it indestructible unless it was originally supercompact.

This theorem is an immediate consequence of the following theorem, where $\coll(\k,\theta)$ is the usual forcing notion which collapses $\theta$ to $\k$. 

\theorem. If $V[G]$ admits a gap below $\k$ and $\k$ is measurable in $V[G]^{\coll(\k,\theta)}$, then $\k$ was $\theta$-supercompact in $V$.

\proof Suppose that $V[G]$ admits a gap below $\k$, that $\k$ remains measurable after the collapse of $\theta$ to $\k$. I must show that $\k$ is $\theta$-supercompact in $V$. Fix an embedding $j:V[G][g]\to M[j(G)][j(g)]$ with critical point $\k$, witnessing that $\k$ remains measurable in $V[G][g]$. Notice that the forcing $G*g$ admits a gap below $\k$. Since moreover $j$ is an ultrapower embedding, $M[j(G)][j(g)]$ is closed under $\k$-sequences in $V[G][g]$. 
Since $\theta$ has been collapsed, it is therefore also closed under $\theta$-sequences. It follows directly now
from the Gap Forcing Theorem of \cite[Ham$\infty$], explained in Remark \FFFGap, that $j\restrict V$ is definable in $V$ and that $M$ is closed under $\theta$-sequences in $V$. Thus, $\k$ is $\theta$-supercompact in $V$, as desired.\qed

\corollary. The following are equivalent:
\points 1. $\k$ is supercompact.\cr
        2. $\k$ is measurable in a forcing extension which admits a gap below $\k$ and in this extension the measurability of $\k$ is indestructible by $\coll(\k,\theta)$ for any $\theta$.\cr 

\proof Certainly $1$ implies $2$ because the Laver preparation (or the Lottery preparation) of a supercompact cardinal $\k$  makes $\k$ indestructible and admits a gap below $\k$. Conversely, $2$ implies $1$ by the previous theorem.\qed

In fact, if the \GCH\ holds, then the result is completely local: $\k$ is $\theta$-supercompact if and only if there is a forcing preparation which admits a gap below $\k$ which makes the measurability of $\k$ indestructible by $\coll(\k,\theta)$. For the forward direction, one can use Theorem \SCLevel. 

While one might suppose from these results that every indestructible strongly compact cardinal is supercompact, this cannot be right because the theorem of Apter and Gitik \cite[AptGit97], which I mentioned earlier, says that it is possible to have a fully indestructible strongly compact cardinal which is also the least measurable cardinal. Such a cardinal could never be supercompact. Beginning with a supercompact cardinal, Apter and Gitik's preparation involves iterated Prikry forcing and consequently does not admit a gap below $\k$. 

The theorem above does show, however, that one cannot hope to make strongly compact non-supercompact cardinals indestructible with forcing that naively resembles the Laver preparation, since all such forcings would admit a gap below $\k$. In particular, one cannot prove that any strongly compact cardinal can be made indestructible by
${\leq}\omega_1$-closed preparatory forcing, or even ${\leq}\omega_1$-strategically closed preparatory forcing, since if
such forcing were prefaced by adding a Cohen real, then the combined forcing would admit a gap. Thus, when it comes to making any strongly compact cardinal fully  indestructible, we evidently need a completely new technique. At the moment, the following questions are open:

\question. Suppose $\k$ is strongly compact. Is there a preparatory forcing to make the strong compactness of $\k$ indestructible by forcing of the form $\add(\k,\d)$? or just by $\add(\k,\k\plus)$?  or of the form $\coll(\k,\d)$? or just $\coll(\k,\k\plus)$?

\medskip
{\parindent=0pt\tenpoint\tightlineskip\sc 
Kobe University, Kobe, Japan, and\par
The City University of New York\par 
\tt hamkins@postbox.csi.cuny.edu\par
http://www.library.csi.cuny.edu/dept/users/hamkins\par}

\quiet\section Bibliography

\input bibliomacros
\tenpoint\tightlineskip

\ref
\author{Arthur W. Apter}
\comment{personal communication}
\key{[Apt96]}

\ref
\author{Arthur W. Apter}
\article{Patterns of Compact Cardinals}
\journal{Annals of Pure and Applied Logic}
\year{1997}
\vol{89}
\no{7}
\page{101-115}
\key{[Apt97]}

\ref
\author{Arthur W. Apter \& Moti Gitik}
\article{The least measurable can be strongly compact and indestructible}
\journal{to appear in the Journal of Symbolic Logic}
\key{[AptGit97]}

\ref
\author{Arthur W. Apter}
\article{Laver indestructibility and the class of compact cardinals}
\journal{Journal of Symbolic Logic}
\year{1998}
\vol{63}
\no{1}
\page{149-157}
\key{[Apt98]}

\ref
\author{James Cummings \& W. Hugh Woodin}
\article{Generalised Prikry forcings}
\comment{(unpublished manuscript)}
\key{[CumWdn]}

\ref
\author{Moti Gitik \& Saharon Shelah}
\article{On certain indestructibility of strong cardinals and a question of 
Hajnal}
\journal{Arch. Math. Logic}
\year{1989}
\vol{28}
\page{35-42}
\key{[GitShl89]}

\ref
\author{Joel David Hamkins}
\article{Lifting and extending measures; fragile measurability}
\year{1994}
\comment{UC Berkeley Dissertation}
\key{[Ham94]}

\ref
\author{Joel David Hamkins}
\article{Canonical seeds and Prikry trees}
\journal{Journal of Symbolic Logic}
\vol{62}
\no{2}
\page{373-396}
\year{1997}
\key{[Ham97]}

\ref
\author{Joel David Hamkins}
\article{Destruction or preservation as you like it}
\journal{Annals of Pure and Applied Logic}
\year{1998}
\vol{91}
\page{191-229}
\key{[Ham98]}

\ref
\author{Joel David Hamkins}
\article{Gap Forcing}
\journal{submitted to the Bulletin of Symbolic Logic}
\comment{(available on the author's web page)}
\key{[Ham$\infty$]}

\ref
\author{K. Kunen, J. Paris}
\article{Boolean extensions and measurable cardinals}
\journal{Annals of Math. Logic}
\vol{2}
\year{1971}
\page{359-377}
\key{[KunPar71]}

\ref
\author{Richard Laver}
\article{Making the supercompactness of $\kappa$ indestructible under 
 $\kappa$-directed closed forcing}
\journal{Israel Journal Math}
\vol{29}
\year{1978}
\page{385-388}
\key{[Lav78]}

\ref
\author{Menachem Magidor}
\article{How large is the first strongly compact cardinal?}
\journal{Annals of Mathematical Logic}
\vol{10}
\year{1976}
\page{33-57}
\key{[Mag76]}

\ref
\author{Telis K. Menas}
\article{On strong compactness and supercompactness}
\journal{Annals of Mathematical Logic}
\vol{7}
\year{1974}
\page{327--359}
\key{[Men74]}

\bye

%% file: fonts.tex
\font\fifteenrm=cmr10 scaled\magstep2 
\font\fifteeni=cmmi10 scaled\magstep2
\font\fifteensy=cmsy10 scaled\magstep2
\font\fifteenbf=cmbx10 scaled\magstep2
\font\fifteentt=cmtt10 scaled\magstep2
\font\fifteenit=cmti10 scaled\magstep2
\font\fifteensl=cmsl10 scaled\magstep2
\font\fifteenam=msam10 scaled\magstep2
\font\fifteenbm=msbm10 scaled\magstep2
\font\fifteenex=cmex10 scaled\magstep2
\font\fifteensc=cmcsc10 scaled\magstep2 
\font\twelverm=cmr10 at 12pt
\font\twelvei=cmmi10 at 12pt
\font\twelvesy=cmsy10 at 12pt
\font\twelvebf=cmbx10 at 12pt
\font\twelvett=cmtt10 at 12pt
\font\twelveit=cmti10 at 12pt
\font\twelvesl=cmsl10 at 12pt
\font\twelveam=msam10 at 12pt
\font\twelvebm=msbm10 at 12pt
\font\twelveex=cmex10 at 12pt
\font\twelvesc=cmcsc10 at 12pt
\font\elevenrm=cmr10 scaled\magstephalf 
\font\eleveni=cmmi10 scaled\magstephalf
\font\elevensy=cmsy10 scaled\magstephalf
\font\elevenbf=cmbx10 scaled\magstephalf
\font\eleventt=cmtt10 scaled\magstephalf
\font\elevenit=cmti10 scaled\magstephalf
\font\elevensl=cmsl10 scaled\magstephalf
\font\elevenam=msam10 scaled\magstephalf
\font\elevenbm=msbm10 scaled\magstephalf
\font\elevenex=cmex10 scaled\magstephalf
\font\elevensc=cmcsc10 scaled\magstephalf
\font\tenrm=cmr10
\font\teni=cmmi10
\font\tensy=cmsy10
\font\tenbf=cmbx10
\font\tentt=cmtt10
\font\tenit=cmti10
\font\tensl=cmsl10
\font\tenam=msam10
\font\tenbm=msbm10
\font\tenex=cmex10
\font\tensc=cmcsc10
\font\ninerm=cmr9
\font\ninei=cmmi9
\font\ninesy=cmsy9
\font\ninebf=cmbx9
\font\ninett=cmtt9
\font\nineit=cmti9
\font\ninesl=cmsl9
\font\nineam=msam9
\font\ninebm=msbm9
\font\nineex=cmex9
\font\ninesc=cmcsc9
\font\eightrm=cmr8
\font\eighti=cmmi8
\font\eightsy=cmsy8
\font\eightbf=cmbx8
\font\eighttt=cmtt8
\font\eightit=cmti8
\font\eightsl=cmsl8
\font\eightam=msam8
\font\eightbm=msbm8
\font\eightex=cmex8
\font\eightsc=cmcsc8
\font\sevenrm=cmr7
\font\seveni=cmmi7
\font\sevensy=cmsy7
\font\sevenbf=cmbx7

\font\sevenam=msam7
\font\sevenbm=msbm7

\font\sixrm=cmr6
\font\sixi=cmmi6
\font\sixsy=cmsy6

\font\sixam=msam6
\font\sixbm=msbm6

\font\fiverm=cmr5
\font\fivei=cmmi5
\font\fivesy=cmsy5

\font\fiveam=msam5
\font\fivebm=msbm5

\font\fourrm=cmr5 at 4pt
\font\fouri=cmmi5 at 4pt
\font\foursy=cmsy5 at 4pt

\font\fouram=msam5 at 4pt
\font\fourbm=msbm5 at 4pt

\skewchar\twelvei='177 \skewchar\eleveni='177\skewchar\teni='177
\skewchar\ninei='177 \skewchar\eighti='177\skewchar\seveni='177 
\skewchar\sixi='177 \skewchar\fivei='177 \skewchar\fouri='177
\skewchar\twelvesy='60 \skewchar\elevensy='60 \skewchar\tensy='60
\skewchar\ninesy='60 \skewchar\eightsy='60 \skewchar\sevensy='60 
\skewchar\sixsy='60 \skewchar\fivesy='60 \skewchar\foursy='60
\newfam\itfam
\newfam\slfam
\newfam\bffam
\newfam\ttfam
\newfam\scfam
\newfam\amfam
\newfam\bmfam
\def\eightbig#1{{\hbox{$\left#1\vbox to 6.5pt{}\voidright $}}}
\def\eightBig#1{{\hbox{$\left#1\vbox to 7.5pt{}\voidright $}}}
\def\eightbigg#1{{\hbox{$\left#1\vbox to 10pt{}\voidright $}}}
\def\eightBigg#1{{\hbox{$\left#1\vbox to 13pt{}\voidright $}}}
\def\ninebig#1{{\hbox{$\left#1\vbox to 7.5pt{}\voidright $}}}
\def\nineBig#1{{\hbox{$\left#1\vbox to 8.5pt{}\voidright $}}}
\def\ninebigg#1{{\hbox{$\left#1\vbox to 11.5pt{}\voidright $}}}
\def\nineBigg#1{{\hbox{$\left#1\vbox to 14.5pt{}\voidright $}}}
\def\tenbig#1{{\hbox{$\left#1\vbox to 8.5pt{}\voidright $}}}
\def\tenBig#1{{\hbox{$\left#1\vbox to 9.5pt{}\voidright $}}}
\def\tenbigg#1{{\hbox{$\left#1\vbox to 12.5pt{}\voidright $}}}
\def\tenBigg#1{{\hbox{$\left#1\vbox to 16pt{}\voidright $}}}
\def\elevenbig#1{{\hbox{$\left#1\vbox to 9pt{}\voidright $}}}
\def\elevenBig#1{{\hbox{$\left#1\vbox to 10.5pt{}\voidright $}}}
\def\elevenbigg#1{{\hbox{$\left#1\vbox to 14pt{}\voidright $}}}
\def\elevenBigg#1{{\hbox{$\left#1\vbox to 17.5pt{}\voidright $}}}
\def\twelvebig#1{{\hbox{$\left#1\vbox to 10pt{}\voidright $}}}
\def\twelveBig#1{{\hbox{$\left#1\vbox to 11pt{}\voidright $}}}
\def\twelvebigg#1{{\hbox{$\left#1\vbox to 15pt{}\voidright $}}}
\def\twelveBigg#1{{\hbox{$\left#1\vbox to 19pt{}\voidright $}}}
\def\fifteenbig#1{{\hbox{$\left#1\vbox to 12pt{}\voidright $}}}
\def\fifteenBig#1{{\hbox{$\left#1\vbox to 13.5pt{}\voidright $}}}
\def\fifteenbigg#1{{\hbox{$\left#1\vbox to 18pt{}\voidright $}}}
\def\fifteenBigg#1{{\hbox{$\left#1\vbox to 23pt{}\voidright $}}}
\def\voidright{\right.\nulldelimiterspace=0pt \mathsurround=0pt }
\def\fifteenpoint{
  \textfont0=\fifteenrm \scriptfont0=\twelverm \scriptscriptfont0=\tenrm
  \def\rm{\fam0 \fifteenrm}%
  \textfont1=\fifteeni \scriptfont1=\twelvei \scriptscriptfont1=\teni
  \textfont2=\fifteensy \scriptfont2=\twelvesy \scriptscriptfont2=\tensy
  \textfont3=\fifteenex \scriptfont3=\fifteenex \scriptscriptfont3=\fifteenex
  \def\it{\fam\itfam\fifteenit}\textfont\itfam=\fifteenit
  \def\sl{\fam\slfam\fifteensl}\textfont\slfam=\fifteensl
  \def\bf{\fam\bffam\fifteenbf}\textfont\bffam=\fifteenbf 
    \scriptfont\bffam=\twelvebf\scriptscriptfont\bffam=\tenbf
  \def\tt{\fam\ttfam\fifteentt}\textfont\ttfam=\fifteentt
  \def\sc{\fam\scfam\fifteensc}\textfont\scfam=\fifteensc
  \def\am{\fam\amfam\fifteenam}\textfont\amfam=\fifteenam
    \scriptfont\amfam=\twelveam\scriptscriptfont\amfam=\tenam
  \def\bm{\fam\bmfam\fifteenbm}\textfont\bmfam=\fifteenbm
    \scriptfont\bmfam=\twelvebm\scriptscriptfont\bmfam=\tenbm
  \baselineskip=21pt \rm
  \let\big=\fifteenbig\let\Big=\fifteenBig\let\bigg=\fifteenbigg
  \let\Bigg=\fifteenBigg}
\def\twelvepoint{
  \textfont0=\twelverm \scriptfont0=\ninerm \scriptscriptfont0=\sevenrm
  \def\rm{\fam0 \twelverm}%
  \textfont1=\twelvei \scriptfont1=\ninei \scriptscriptfont1=\seveni
  \textfont2=\twelvesy \scriptfont2=\ninesy \scriptscriptfont2=\sevensy
  \textfont3=\twelveex \scriptfont3=\twelveex \scriptscriptfont3=\twelveex
  \def\it{\fam\itfam\twelveit}\textfont\itfam=\twelveit
  \def\sl{\fam\slfam\twelvesl}\textfont\slfam=\twelvesl
  \def\bf{\fam\bffam\twelvebf}\textfont\bffam=\twelvebf 
    \scriptfont\bffam=\ninebf\scriptscriptfont\bffam=\sevenbf
  \def\tt{\fam\ttfam\twelvett}\textfont\ttfam=\twelvett
  \def\sc{\fam\scfam\twelvesc}\textfont\scfam=\twelvesc
  \def\am{\fam\amfam\twelveam}\textfont\amfam=\twelveam
    \scriptfont\amfam=\nineam\scriptscriptfont\amfam=\sevenam
  \def\bm{\fam\bmfam\twelvebm}\textfont\bmfam=\twelvebm
    \scriptfont\bmfam=\ninebm\scriptscriptfont\bmfam=\sevenbm
  \baselineskip=17.8pt \rm 
  \def\looselineskip{\baselineskip=18.5pt plus 1.8pt}%
  \def\tightlineskip{\baselineskip=16.5pt}%
  \def\verytightlineskip{\baselineskip=15pt}%
  \let\big=\twelvebig\let\Big=\twelveBig\let\bigg=\twelvebigg
  \let\Bigg=\twelveBigg  }
\def\elevenpoint{
  \textfont0=\elevenrm \scriptfont0=\ninerm \scriptscriptfont0=\sixrm
  \def\rm{\fam0 \elevenrm}%
  \textfont1=\eleveni \scriptfont1=\ninei \scriptscriptfont1=\sixi
  \textfont2=\elevensy \scriptfont2=\ninesy \scriptfont2=\sixsy 
  \textfont3=\elevenex \scriptfont3=\elevenex \scriptfont3=\elevenex
  \def\it{\fam\itfam\elevenit}\textfont\itfam=\elevenit
  \def\sl{\fam\slfam\elevensl}\textfont\slfam=\elevensl
  \def\bf{\fam\bffam\elevenbf}\textfont\bffam=\elevenbf
  \def\tt{\fam\ttfam\eleventt}\textfont\ttfam=\eleventt
  \def\sc{\fam\scfam\elevensc}\textfont\scfam=\elevensc
  \def\am{\fam\amfam\elevenam}\textfont\amfam=\elevenam
    \scriptfont\amfam=\nineam\scriptscriptfont\amfam=\sixam
  \def\bm{\fam\bmfam\elevenbm}\textfont\bmfam=\elevenbm
    \scriptfont\bmfam=\ninebm\scriptscriptfont\bmfam=\sixbm
  \baselineskip=15.1pt \rm
  \def\looselineskip{\baselineskip=16pt plus 1.5pt}%
  \def\tightlineskip{\baselineskip=14pt}%
  \def\verytightlineskip{\baselineskip=13pt}%
  \let\big=\elevenbig\let\Big=\elevenBig\let\bigg=\elevenbigg
  \let\Bigg=\elevenBigg  }
\def\tenpoint{
  \textfont0=\tenrm \scriptfont0=\eightrm \scriptscriptfont0=\fiverm
  \def\rm{\fam0 \tenrm}%
  \textfont1=\teni \scriptfont1=\eighti \scriptscriptfont1=\fivei
  \textfont2=\tensy \scriptfont2=\eightsy \scriptfont2=\fivesy 
  \textfont3=\tenex \scriptfont3=\tenex \scriptfont3=\tenex
  \def\it{\fam\itfam\tenit}\textfont\itfam=\tenit
  \def\sl{\fam\slfam\tensl}\textfont\slfam=\tensl
  \def\bf{\fam\bffam\tenbf}\textfont\bffam=\tenbf
  \def\tt{\fam\ttfam\tentt}\textfont\ttfam=\tentt
  \def\sc{\fam\scfam\tensc}\textfont\scfam=\tensc
  \def\am{\fam\amfam\tenam}\textfont\amfam=\tenam
    \scriptfont\amfam=\eightam \scriptscriptfont\amfam=\fiveam
  \def\bm{\fam\bmfam\tenbm}\textfont\bmfam=\tenbm
    \scriptfont\bmfam=\eightbm \scriptscriptfont\bmfam=\fivebm
  \baselineskip=14pt \rm
  \def\looselineskip{\baselineskip=14.8pt plus1.5pt}
  \def\tightlineskip{\baselineskip=13.6pt}%
  \def\verytightlineskip{\baselineskip=13pt}%
  \let\big=\tenbig\let\Big=\tenBig\let\bigg=\tenbigg\let\Bigg=\tenBigg  }
\def\ninepoint{
  \textfont0=\ninerm \scriptfont0=\sevenrm \scriptscriptfont0=\fourrm
  \def\rm{\fam0 \ninerm}%
  \textfont1=\ninei \scriptfont1=\seveni \scriptscriptfont1=\fouri
  \textfont2=\ninesy \scriptfont2=\sevensy \scriptfont2=\foursy 
  \textfont3=\nineex \scriptfont3=\nineex \scriptfont3=\nineex
  \def\it{\fam\itfam\nineit}\textfont\itfam=\nineit
  \def\sl{\fam\slfam\ninesl}\textfont\slfam=\ninesl
  \def\bf{\fam\bffam\ninebf}\textfont\bffam=\ninebf
  \def\tt{\fam\ttfam\ninett}\textfont\ttfam=\ninett
  \def\sc{\fam\scfam\ninesc}\textfont\scfam=\ninesc
  \def\am{\fam\amfam\nineam}\textfont\amfam=\nineam
    \scriptfont\amfam=\nineam\scriptscriptfont\amfam=\fouram
  \def\bm{\fam\bmfam\ninebm}\textfont\bmfam=\ninebm
    \scriptfont\bmfam=\ninebm\scriptscriptfont\bmfam=\fourbm
  \baselineskip=12.6pt \rm
  \let\big=\ninebig\let\Big=\nineBig\let\bigg=\ninebigg
  \let\Bigg=\nineBigg  }
\def\eightpoint{
  \textfont0=\eightrm \scriptfont0=\fiverm \scriptscriptfont0=\fourrm
  \def\rm{\fam0 \eightrm}%
  \textfont1=\eighti \scriptfont1=\fivei \scriptscriptfont1=\fouri
  \textfont2=\eightsy \scriptfont2=\fivesy \scriptfont2=\foursy 
  \textfont3=\eightex \scriptfont3=\eightex \scriptfont3=\eightex
  \def\it{\fam\itfam\eightit}\textfont\itfam=\eightit
  \def\sl{\fam\slfam\eightsl}\textfont\slfam=\eightsl
  \def\bf{\fam\bffam\eightbf}\textfont\bffam=\eightbf
  \def\tt{\fam\ttfam\eighttt}\textfont\ttfam=\eighttt
  \def\sc{\fam\scfam\eightsc}\textfont\scfam=\eightsc
  \def\am{\fam\amfam\eightam}\textfont\amfam=\eightam
    \scriptfont\amfam=\eightam\scriptscriptfont\amfam=\fouram
  \def\bm{\fam\bmfam\eightbm}\textfont\bmfam=\eightbm
    \scriptfont\bmfam=\eightbm\scriptscriptfont\bmfam=\fourbm
  \baselineskip=11.2pt \rm
  \let\big=\eightbig\let\Big=\eightBig\let\bigg=\eightbigg
  \let\Bigg=\eightBigg  }

%% file: articlemacros.tex
\twelvepoint
\nopagenumbers
\hsize=6in\vsize=8.8in

\parskip=1pt plus 1pt

\newif\ifSpecialhead\Specialheadfalse
\newbox\specialheadbox

\def\specialhead #1\par{\Specialheadtrue\setbox\specialheadbox=\hbox{#1}}
\headline={{\ifSpecialhead\box\specialheadbox\global\Specialheadfalse\else
     \ifnum\pageno<0{\hfill\quad{\twelvebf\folio}}%
     \else\ifnum\pageno<2\hfill
     \else\hfill\twelvepoint\sc\firstmark\quad{\twelvebf\folio}\fi\fi\fi}}

\def\title#1\par{\bigskip{\def\cr{\par\center}\center\fifteenbf #1\par}\medskip}
\def\subtitle#1\par{\centerline{\fifteenrm #1}\medskip}
\def\author#1\par{\medskip{\def\cr{\par\center\twelvesc}\fifteensc\center#1\par}}
\def\center#1\par{\hfil #1\hfil\par}
\def\abstract.#1\par{\message{Abstract.}%
                    \medskip{\narrower\narrower\tenpoint\tightlineskip
                        \noindent{\bf Abstract.}#1\par}\medskip\noindent}
\def\bigabstract.#1\par{\message{Abstract.}%
                         \medskip{\narrower\narrower\tightlineskip
                         \noindent{\bf Abstract. }#1\par}\medskip\noindent}
\def\acknowledgement#1\par{\footnote{}{#1}}
\def\sectionskip{\Goodbreak\vskip 25pt plus 15pt minus 5pt}
\def\secnumber{\ifquiet
               \else\ifNoSections
                    \else\sectionsymbol\the\secno\quad\fi\fi}
\def\section#1\par{ \NoSectionsfalse\par\sectionskip\proofdepth=0\claimno=0
 \ifquiet\else\advance\secno by1\fi\toks0={#1}
 \immediate\write16{\ifquiet\else Section \the\secno\space\fi
                    \the\toks0}%
 \mark{\secnumber #1}%
 {\fifteenpoint\bf\noindent\secnumber #1}\nobreak\bigskip\quietoff
 \nobreak\noindent}
\def\quiet{\quiettrue}

\def\quietoff{\ifQUIET\else\quietfalse\fi}
\newif\ifquiet
\newif\ifQUIET
\newif\ifNoSections
\newcount\claimtype
\newcount\secno
\newcount\claimno
\newcount\subclaimno
\newcount\subsubclaimno
\newcount\subsubsubclaimno
\newcount\proofdepth
\def\subclaimnumber{\ifquiet\else\ifcase\subclaimno\or A\or B\or C\or D\or E\or
     F\or G\or H\or I\or J\or K\or L\or M\or N\or O\or P\fi\fi}
\def\subsubclaimnumber{\ifquiet\else\ifcase\subsubclaimno\or i\or ii\or iii\or 
   iv\or v\or vi\or vii\or viii\or ix\or x\or xi\or xii\or xiii\or xiv\fi\fi}
\def\subsubsubclaimnumber{\ifquiet\else\ifcase\subsubsubclaimno\or a\or b\or 
   c\or d\or e\or f\or g\or \or h\or i\or j\or k\or l\or m\or n\or o\fi\fi}
\def\claimtag{\ifquiet\else
  \ifNoSections
    \ifcase\proofdepth\the\claimno%
    \or\the\claimno.\subclaimnumber
    \or\the\claimno.\subclaimnumber.\subsubclaimnumber
    \or\the\claimno.\subclaimnumber.\subsubclaimnumber
                                                .\subsubsubclaimnumber\fi
  \else
    \ifcase\proofdepth\the\secno.\the\claimno
    \or\the\secno.\the\claimno.\subclaimnumber
    \or\the\secno.\the\claimno.\subclaimnumber.\subsubclaimnumber
    \or\the\secno.\the\claimno.\subclaimnumber.\subsubclaimnumber
                                                .\subsubsubclaimnumber\fi\fi\fi}
\secno=0\claimno=0\proofdepth=0\subclaimno=0\subsubclaimno=0\subsubsubclaimno=0
\NoSectionstrue
\newbox\qedbox
\def\claimname{\ifcase\claimtype Theorem\or Lemma\or Claim\or Corollary\or
               Question\or Definition\or Remark\or Conjecture\fi}
\def\preclaimskip{\removelastskip
    \ifcase\claimtype\goodbreak\vskip 8pt plus 4pt minus 2pt
                  \or\goodbreak\vskip 6pt plus 4pt minus 1pt
                  \or\goodbreak\vskip 5pt plus 4pt minus 1pt
                  \or\goodbreak\vskip 8pt plus 4pt minus 2pt
                  \or\vskip 7pt plus 4pt minus 2pt
                  \or\vskip 7pt plus 4pt minus 2pt
                  \or\vskip 7pt plus 4pt minus 2pt
                  \or\goodbreak\vskip 8pt plus 4pt minus 2pt\fi}
\def\postclaimskip{\ifcase\claimtype         \vskip 4pt plus 2pt minus 2pt
                                          \or\vskip 3pt plus 2pt minus 2pt
                                          \or\vskip 2pt plus 2pt minus 1pt
                                          \or\vskip 4pt plus 2pt minus 2pt
                                          \or\vskip 1pt plus 2pt 
                                          \or\vskip 4pt plus 4pt 
                                          \or\vskip 3pt plus 2pt
                                          \or\vskip 4pt plus 2pt minus 2pt\fi}
\def\claimfont{\ifcase\claimtype
                  \sl\or\sl\or\sl\or\sl\or\sl\or\rm\or\rm\or\sl\fi}
\def\advancetag{\ifcase\proofdepth\advance\claimno by1
                               \or\advance\subclaimno by1
                               \or\advance\subsubclaimno by1
                               \or\advance\subsubsubclaimno by1\fi}
\def\sayclaim#1.#2 #3\par{\ifquiet\else\advancetag\fi
    \preclaimskip\setbox1=\hbox{#1}\setbox2=\hbox{#2}%
    \toks0={#1 }
    \immediate\write16{\ifdim\wd1>0pt\the\toks0
                       \else\claimname\space\fi \claimtag.}%
    \vbox{\noindent
    {\bf\ifdim\wd1=0pt \claimname\else #1\fi\ifquiet.\else\ \claimtag{\ifNoSections.\fi}\fi}%
    \enspace{\ifdim\wd2>0pt\sc #2\enspace\fi}%
    {\claimfont #3\par}}\postclaimskip\quietoff}
\def\theorem{\claimtype=0\sayclaim}
\def\lemma{\claimtype=1\sayclaim}

\def\corollary{\claimtype=3\sayclaim}
\def\question{\claimtype=4\sayclaim}

\def\remark{\claimtype=6\sayclaim}

\def\point#1. #2\par{\item{\rm #1.}#2\par}
\def\points#1\cr\par{\medskip\vbox{\let\cr=\point\point#1\par}\par}
\def\df{\it}
\def\prooffont{}
\def\proofsize{}
\def\proofindent{}
\def\proofskip{\badbreak\ifcase\claimtype    \vskip 3pt plus 2pt minus 2pt
                                          \or\vskip 2pt plus 2pt minus 2pt
                                          \or\vskip 1pt plus 2pt minus 1pt
                                          \or\vskip 3pt plus 2pt minus 2pt
                                          \or\vskip 1pt plus 2pt 
                                          \or\vskip 2pt plus 4pt 
                                          \or\vskip 1pt plus 2pt
                                          \or\vskip 3pt plus 2pt minus 2pt\fi}

\def\Goodbreak{\vskip0pt plus.5in\penalty-1000\vskip0pt plus-.5in}
\def\goodbreak{\penalty-500}
\def\badbreak{\penalty500}
\def\Badbreak{\penalty1000}
\def\proof{\message{proof}\removelastskip\Badbreak\proofskip\begingroup
  \advance\proofdepth by1
  \setbox\qedbox=\hbox{\halmos\raise2pt\hbox{\fiverm\claimname}}%
  \prooffont\proofsize\proofindent\noindent{\bf Proof: }}
\def\proofof#1:{\message{proof}\removelastskip\Badbreak\proofskip\begingroup
  \advance\proofdepth by1
  \setbox\qedbox=\hbox{\halmos\raise2pt\hbox{\fiverm#1}}%
  \prooffont\proofsize\proofindent\noindent{\bf Proof of #1: }}
\def\cite[#1]{[{\tenrm{#1}}]\message{[#1]}}
\edef\ref#1{\expandafter\global\expandafter\edef#1{\noexpand\claimtag}}
\newwrite\notes
\openout\notes=\jobname.notes
\long\def\unexpandedwrite#1#2{\def\finwrite{\write#1}%
   {\aftergroup\finwrite\aftergroup{\sanitize#2\endsanity}}}
\def\sanitize{\futurelet\next\sanswitch}
\let\stoken=\space
\def\sanswitch{\ifx\next\endsanity
  \else\ifcat\noexpand\next\stoken\aftergroup\space\let\next=\eat
   \else\ifcat\noexpand\next\bgroup\aftergroup{\let\next=\eat
    \else\ifcat\noexpand\next\egroup\aftergroup}\let\next=\eat
     \else\let\next=\copytoken\fi\fi\fi\fi \next}
\def\eat{\afterassignment\sanitize \let\next= }
\long\def\copytoken#1{\ifcat\noexpand#1\relax\aftergroup\noexpand
  \else\ifcat\noexpand#1\noexpand~\aftergroup\noexpand\fi\fi
  \aftergroup#1\sanitize}
\def\endsanity\endsanity{}

\def\note#1#2{\hbox to2in{\strut#1\quad\dotfill\quad#2}}
\def\boxit#1{\setbox4=\hbox{\kern1pt#1\kern1pt}
  \hbox{\vrule\vbox{\hrule\kern1pt\box4\kern1pt\hrule}\vrule}}
\def\halmos{\hbox{\am\char'3}} 
\def\qed#1\par{\message{.                                }\setbox1=\hbox{#1}%
  \ifdim\wd1>0pt\setbox\qedbox=\hbox{\halmos\raise2pt\hbox{\fiverm #1}}\fi
  \kern5pt\lower 2pt\hbox{\box\qedbox}\proofskip\goodbreak\endgroup}

%% file: generalmathmacros.tex
\def\sectionsymbol{\S}
\def\k{\kappa}
\def\g{\gamma}
\def\a{\alpha}
\def\b{\beta}
\def\d{\delta}
\def\s{\sigma}
\def\t{\tau}
\def\l{\lambda}

\def\I1{\mathop{\hbox{\sc i}_1}}
\def\w{\omega}
\def\P{{\mathchoice{\hbox{\bm P}}{\hbox{\bm P}}
         {\hbox{\tenbm P}}{\hbox{\sevenbm P}}}}
\def\Q{{\mathchoice{\hbox{\bm Q}}{\hbox{\bm Q}}
         {\hbox{\tenbm Q}}{\hbox{\sevenbm Q}}}}
\def\R{{\mathchoice{\hbox{\bm R}}{\hbox{\bm R}}
         {\hbox{\tenbm R}}{\hbox{\sevenbm R}}}}

\def\F{{\mathchoice{\hbox{\bm F}}{\hbox{\bm F}}
         {\hbox{\tenbm F}}{\hbox{\sevenbm F}}}}
\def\X{{\mathchoice{\hbox{\bm X}}{\hbox{\bm X}}
         {\hbox{\tenbm X}}{\hbox{\sevenbm X}}}}
\def\card#1{\left|#1\right|}

\def\dom{\mathop{\rm dom}\nolimits}

\def\coll{\mathop{\rm coll}}

\def\id{\mathop{\hbox{\tenrm id}}}

\def\elesub{\prec}

\def\unifto{\buildrel\lower 7pt\hbox{$\to$}\over\to}

\def\iso{\cong}

\def\cof{\mathop{\rm cof}\nolimits}
\def\cp{\mathop{\rm cp}\nolimits}
\def\ran{\mathop{\rm ran}\nolimits}
\def\from{\mathbin{\vbox{\baselineskip=3pt\lineskiplimit=0pt
                         \hbox{.}\hbox{.}\hbox{.}}}}
\def\ORD{\hbox{\sc ord}}

\def\HOD{\hbox{\sc hod}}

\def\GCH{\hbox{\sc gch}}

\def\inacc{\hbox{\sc inacc}}

\def\plus{^{\scriptscriptstyle +}}
\def\plusplus{^{\scriptscriptstyle ++}}

\def\in{\mathrel{\mathchoice{\raise 
1pt\hbox{$\scriptstyle\cal\char'62$}}
         {\raise 1pt\hbox{$\scriptstyle\cal\char'62$}}
         {\raise .5pt\hbox{$\scriptscriptstyle\cal\char'62$}}
         {\hbox{$\scriptscriptstyle\cal\char'62$}}}\penalty700{}}
\def\ni{\mathrel{\mathchoice{\raise 1pt\hbox{$\scriptstyle\cal\char'63$}}
                   {\raise 1pt\hbox{$\scriptstyle\cal\char'63$}}
                   {\raise .5pt\hbox{$\scriptscriptstyle\cal\char'63$}}
                   {\hbox{$\scriptscriptstyle\cal\char'63$}}}\penalty700}
\def\of{\mathrel{\mathchoice{\raise 1pt\hbox{$\scriptstyle\subseteq$}}
                   {\raise 1pt\hbox{$\scriptstyle\subseteq$}}
                   {\raise .5pt\hbox{$\scriptscriptstyle\subseteq$}}
                   {\hbox{$\scriptscriptstyle\subseteq$}}}}
\def\fo{\mathrel{\mathchoice{\raise 1pt\hbox{$\scriptstyle\supseteq$}}
                   {\raise 1pt\hbox{$\scriptstyle\supseteq$}}
                   {\raise .5pt\hbox{$\scriptscriptstyle\supseteq$}}
                   {\hbox{$\scriptscriptstyle\supseteq$}}}}
\def\notin{\mathrel{\mathchoice
  {\raise 1pt\hbox{\rlap{$\scriptstyle\;|$}$\scriptstyle\cal\char'62$}}
  {\raise 1pt\hbox{\rlap{$\scriptstyle\kern2pt 
          |$}$\scriptstyle\cal\char'62$}}
  {\raise .5pt\hbox{\rlap{$\scriptscriptstyle\, |$}$\scriptscriptstyle
      \cal\char'62$}}
  {\hbox{\rlap{$\scriptscriptstyle\, |$}$\scriptscriptstyle
     \cal\char'62$}}}%
  \penalty700}

\def\and{\mathrel{\kern1pt\&\kern1pt}}
\def\iff{\mathrel{\leftrightarrow}}

\def\union{\cup}

\def\compose{\circ}
\def\intersect{\cap}

\def\setminus{\mathbin{\hbox{\bm\char'162}}}

\def\add{\mathop{\rm Add}\nolimits}

\def\cross{\times}

\def\lte{\mathrel{\scriptstyle\leq}}

\def\[#1]{\left[\vphantom{\bigm|}#1\right]}
\def\<#1>{\langle\,#1\,\rangle}

\def\image{\mathbin{\hbox{\tt\char'42}}}
\def\restrict{\mathbin{\mathchoice{\hbox{\am\char'26}}{\hbox{\am\char'26}}{\hbox{\eightam\char'26}}{\hbox{\sixam\char'26}}}}
\def\force{\mathbin{\hbox{\am\char'15}}}

\def\beth{\mathord{\hbox{\bm\char'151}}}

\def\st{\mid}
\def\seq<#1>{{\def\st{\mid\penalty650}\left<\,#1\,\right>}}
\def\sing#1{\{#1\}}
\def\set#1{\{\,#1\,\}}

\def\th{{\hbox{\fiverm th}}}

\def\forces{\force}
\def\lttheta{{\raise 1pt\hbox{$\scriptstyle<$}\theta}}

\def\I1{\mathop{\hbox{\sc i}_1}}
\def\ltk{{{\scriptstyle<}\k}}

\def\ltg{{{\scriptstyle<}\g}}
\def\ltd{{{\scriptstyle<}\d}}
\def\ltb{{{\scriptstyle<}\b}}
\def\lteb{{{\scriptstyle\leq}\b}}
\def\lted{{{\scriptstyle\leq}\d}}
\def\lteg{{{\scriptstyle\leq}\g}}
\def\ltet{{{\scriptstyle\leq}\theta}}
\def\ltel{{{\scriptstyle\leq}\l}}

\def\hdot{\dot h}

\def\one{1\hskip-3pt {\rm l}}
\def\Ptail{\P_{\fiverm \!tail}}

\def\Xdot{\dot X}
\def\Ydot{\dot Y}
\def\ltg{{\scriptscriptstyle<}\g}

\def\Gtail{G_{\fiverm tail}}
\def\Ftail{\F_{\!\fiverm tail}}
\def\ftail{f_{\fiverm tail}}

\def\ltek{{\scriptstyle\leq}\k}

\def\Vbar{{\overline V}}
\def\Mbar{{\overline M}}
\def\Dbar{{\overline D}}
\def\Cbar{{\overline C}}
\def\Qtilde{{\widetilde \Q}}
\def\Gtilde{{\widetilde G}}
\def\Ptilde{{\widetilde\P}}
\def\Diamond{\diamondsuit}

%% file: arrowmacros.tex
\font\arrow=line10 scaled \magstep1
\def\makeline#1.{\hbox{\arrow\char#1}}
\def\makearrow#1.#2.{\hbox{\arrow\char#1\llap{\char#2}}}
\def\definelinesandarrows#1.#2.#3.#4.#5.{
   \expandafter\edef\csname#4line\endcsname{\makeline#1.}
   \expandafter\edef\csname#4arrow\endcsname{\makearrow#1.#2.}
   \expandafter\edef\csname#5line\endcsname{\makeline#1.}
   \expandafter\edef\csname#5arrow\endcsname{\makearrow#1.#3.}}
\definelinesandarrows 0.18.9.ne.sw.
\definelinesandarrows 1.21.11.nnne.sssw.
\definelinesandarrows 2.14.13.nnnne.ssssw.
\definelinesandarrows 3.23.15.nnnnne.sssssw.
\definelinesandarrows 4.23.15.nnnnnne.ssssssw.
\definelinesandarrows 10.30.29.nne.ssw.
\definelinesandarrows 16.49.41.neeeeee.swwwwww.
\definelinesandarrows 17.51.43.neeee.swwww.
\definelinesandarrows 19.55.47.nehuh.swhuh.
\definelinesandarrows 24.58.41.neeeeeee.swwwwwww.
\definelinesandarrows 26.62.9.neee.swww.
\definelinesandarrows 33.49.25.neeeee.swwwww.
\definelinesandarrows 35.62.61.nee.sww.
\definelinesandarrows 64.82.73.se.nw.
\definelinesandarrows 65.85.75.ssse.nnnw.
\definelinesandarrows 66.78.77.sssse.nnnnw.
\definelinesandarrows 67.87.79.ssssse.nnnnnw.
\definelinesandarrows 68.87.79.sssssse.nnnnnnw.
\definelinesandarrows 74.94.93.sse.nnw.
\definelinesandarrows 80.113.105.seeeeee.nwwwwww.
\definelinesandarrows 81.115.107.seeee.nwwww.
\definelinesandarrows 99.126.125.see.nww.
\def\sejoin#1#2{\setbox1=\hbox{#1}\setbox2=\hbox{#2}%
  \hbox{\vbox{\hbox{\copy1\kern\wd2}\nointerlineskip
              \hbox{\kern\wd1\box2}}}}
\def\nejoin#1#2{\setbox1=\hbox{#1}\setbox2=\hbox{#2}%
  \hbox{\vbox{\hbox{\kern\wd1\copy2}\nointerlineskip\hbox{\copy1\kern\wd2}}}}
\newdimen\hnudge
\newdimen\vnudge
\newdimen\hnudgedefault
\newdimen\vnudgedefault

\def\SEdefaultnudge{\hnudge=-16pt\vnudge=20pt}
\def\Edefaultnudge{\hnudge=-25pt\vnudge=6pt}
\def\Sdefaultnudge{\hnudge=-8pt\vnudge=20pt}
\def\longEdefaultnudge{\hnudge=-5pt\vnudge=6pt}
\def\nudgeright#1pt{\advance\hnudge by#1pt}
\def\nudgeleft#1pt{\advance\hnudge by-#1pt}
\def\nudgeup#1pt{\advance\vnudge by#1pt}
\def\nudgedown#1pt{\advance\vnudge by-#1pt}
\def\label#1{\smash{\llap{\kern\hnudge
                   \raise\vnudge\rlap{$\scriptstyle#1$}\hfill}}}

\def\SEarrow{\SEdefaultnudge
             \sejoin\seeline{\sejoin\seeline{\sejoin\seeline\seearrow}}}

\def\Sarrow{\Sdefaultnudge\setbox1=\hbox{\SEarrow}
           \hbox{\hskip 10pt\vrule height\ht1\hbox{\arrow\char'77}}}
\def\Earrow{\Edefaultnudge\setbox1=\hbox{\SEarrow}
 \hbox{\raise 2pt\hbox{\vrule height-.4pt depth.8ptwidth\wd1\kern2pt
       \llap{\arrow\char'55}}}}
\def\longEarrow{\longEdefaultnudge\setbox1=\hbox{\SEarrow}
      \rlap{\hskip-1.25\wd1\raise 2pt
            \hbox{\vrule height-.4pt depth.8ptwidth2.5\wd1\kern2pt
            \llap{\arrow\char'55}}}}
\def\trianglediagram#1#2#3#4#5#6{%
    {\def\normalbaselines{\baselineskip0pt\lineskip8pt\lineskiplimit0pt}%
       \matrix{#1& &\cr
               \Sarrow\label{#2}&\SEarrow\label{#3}&\cr
               #4&\Earrow\label{#5}&#6\cr}}}

%% file: bibliomacros.tex
\nopagenumbers
\parindent=0pt
\newbox\Article
\newbox\Journal
\newbox\Author
\newbox\Vol
\newbox\No
\newbox\Year
\newbox\Page
\newbox\Book
\newbox\Publisher
\newbox\Pubaddr
\newbox\Key
\newbox\Editor
\newbox\Comment
\newbox\Note
\def\entry#1#2\par{\item{#1\quad}\hskip-1.1em#2\par}
\def\article#1{\setbox\Article=\hbox{\sl #1, }}
\def\journal#1{\setbox\Journal=\hbox{\rm #1 }}
\def\author#1{\setbox\Author=\hbox{\sc #1, }}
\def\vol#1{\setbox\Vol=\hbox{\bf #1 }}
\def\no#1{\setbox\No=\hbox{no. #1 }}
\def\year#1{\setbox\Year=\hbox{\rm({\oldstyle #1}) }}
\def\page#1{\setbox\Page=\hbox{\rm p. #1 }}
\def\book#1{\setbox\Book=\hbox{\it #1, }}
\def\publisher#1{\setbox\Publisher=\hbox{\rm #1, }}
\def\pubaddr#1{\setbox\Pubaddr=\hbox{\rm #1, }}
\def\key#1{\setbox\Key=\hbox{#1}}
\def\editor#1{\setbox\Editor=\hbox{\rm(#1, Ed.) }}
\def\comment#1{\setbox\Comment=\hbox{\rm #1}}
\def\note#1{\setbox\Note=\hbox{\rm #1 }}
\def\ref#1\par{\smallskip{#1
  \entry{\ifhbox\Key\unhbox\Key\else[\ ]\fi}%
  \unhbox\Author\unhbox\Note
  \ifhbox\Book \unhbox\Book\unhbox\Publisher\unhbox\Pubaddr
               \unhbox\Editor\unhbox\Page\unhbox\Year\unhbox\Comment
  \else \unhbox\Article\unhbox\Journal\unhbox\Vol\unhbox\No\unhbox\Editor
        \unhbox\Page\unhbox\Year\unhbox\Comment\fi\par}}